\documentclass
[
11pt,
]
{amsart}
\usepackage{
amssymb
}
\newcommand{\no}[1]{#1}

\renewcommand{\no}[1]{}\newcommand{\wbar}{\bar} \newcommand{\wtilde}{\tilde} \newcommand{\upDelta}{\Delta} 

\no{\usepackage{times}\usepackage[
subscriptcorrection, slantedGreek, 
nofontinfo]{mtpro}
\renewcommand{\Delta}{\upDelta}
}

\newtheorem{theorem}{Theorem}
\newtheorem{proposition}{Proposition}
\newtheorem{lemma}{Lemma}
\newtheorem{definition}{Definition}
\newtheorem{corollary}{Corollary}

\DeclareMathOperator{\supp}{supp}
\DeclareMathOperator{\Ker}{Ker}
\newcommand{\eps}{\varepsilon}
\newcommand{\R}{{\bf R}}
\newcommand{\Id}{\mbox{Id}}
\renewcommand{\r}[1]{(\ref{#1})}
\newcommand{\PDO}{$\Psi$DO}
\newcommand{\be}[1]{\begin{equation}\label{#1}}
\newcommand{\ee}{\end{equation}}

\DeclareMathOperator{\n}{neigh}
\renewcommand{\d}{\mathrm{d}}
\renewcommand{\i}{\mathrm{i}}

\newcommand{\bo}{\partial M}

\newcommand{\Mint}{M^\text{\rm int}}

\title[Integral geometry of tensor fields]{Integral geometry of tensor fields on a class of non-simple Riemannian manifolds}

\author[P. Stefanov]{Plamen Stefanov}
\address{Department of Mathematics, Purdue University, West Lafayette, IN 47907}
\thanks{First author partly supported by NSF Grant DMS-0196440}
\author[G. Uhlmann]{Gunther Uhlmann}
\address{Department of Mathematics, University of Washington, Seattle, WA 98195}
\thanks{Second author partly supported by NSF and a John Simon Guggenheim fellowship}

\begin{document}

\begin{abstract} 
We study the geodesic X-ray transform $I_\Gamma$ of tensor fields on a compact Riemannian manifold $M$ with non-necessarily convex boundary and with possible conjugate points. We assume that $I_\Gamma$ is known for geodesics belonging to an open set $\Gamma$ with endpoints on the boundary. 
We prove generic s-injectivity and a stability estimate under some topological assumptions and under the condition that for any $(x,\xi)\in T^*M$, there is a geodesic in $\Gamma$ through $x$ normal to $\xi$ without conjugate points. 
\end{abstract}

\maketitle

\section{Introduction and statement of the main results} 
Let $(M,\partial M)$ be a smooth compact  manifold with boundary, and let $g\in C^k(M)$ be a Riemannian metric on it. We can always assume that $(M,\partial M)$ is equipped with a real analytic atlas, while $\partial M$ and $g$ may or may not be analytic. 
We define the geodesic X-ray transform $I$ of symmetric 2-tensor fields by
\be{I_G}
I f(\gamma) = \int_0^{l_\gamma} \langle f(\gamma(t)), \dot \gamma^2(t) \rangle \,\d t, 
\ee
where $[0,l_\gamma]\ni t\mapsto \gamma$ is any  geodesic with endpoints on $\bo$ parameterized by its arc-length. 
Above, $\langle f, \theta^2\rangle$ is the action of $f$ on the vector $\theta$, that in local coordinates is given by $f_{ij}  \theta^i \theta^j$. The purpose of this work is to study the injectivity, up to potential fields,  and stability estimates for $I$ restricted to certain subsets $\Gamma$ (that we call $I_\Gamma$), and for manifolds with possible conjugate points. We require however that the geodesics in $\Gamma$ do not have conjugate points. We also require that $\Gamma$ is an open sets of geodesics such that the collection of their conormal bundles covers $T^*M$. This guarantees that $I_\Gamma$ resolves the singularities. The main results are injectivity up to a potential field and stability for generic metrics, and in particular for real analytic ones.

We are motivated here by the boundary rigidity problem: to recover $g$, up to an isometry leaving $\partial M$ fixed, from knowledge of the boundary distance function $\rho(x,y)$ for a subset of pairs $(x,y)\in \partial M\times \partial M$, see e.g., \cite{M,Sh, CDS, SU-rig,PU}.  In presence of conjugate points, one should study instead the lens rigidity problem: a recovery of $g$ from its scattering relation restricted to a subset. Then $I_\Gamma$ is the linearization of those problems for an appropriate $\Gamma$. Since we want to trace the dependence of $I_\Gamma$ on perturbations of the metric, it is more convenient to work with open $\Gamma$'s that have dimension larger than $n$, if $n\ge3$, making the linear inverse problem formally overdetermined. One can use the same method to study restrictions of $I$ on $n$ dimensional subvarieties but this is behind the scope of this work. 

Any symmetric 2-tensor field $f$ can be written as an orthogonal sum of a \textit{solenoidal} part $f^s$ and a \textit{potential} one $dv$, where $v=0$ on $\bo$, and $d$ stands for the symmetric differential of the 1-form $v$, see Section~\ref{sec_prel}. Then $I(dv)(\gamma)=0$ for any geodesic $\gamma$ with endpoints on $\bo$. We say that $I_\Gamma$ is \textit{s-injective}, if $I_\Gamma f=0$ implies $f=dv$ with $v=0$ on $\bo$, or, equivalently, $f=f^s$. 
This problem has been studied before for \textit{simple} manifolds with boundary, i.e., under the assumption that $\bo$ is strictly convex, and there are no conjugate points in $M$ (then $M$ is diffeomorphic to a ball).   The book \cite{Sh} contains the main results
up to 1994 on the integral geometry problem considered in this paper. Some recent results include \cite{Sh-sib}, \cite{Ch}, \cite{SU-Duke}, \cite{D}, \cite{Pe}, \cite{SSU}, \cite{SU}. In the two dimensional case, following the method used in \cite{PU} to solve the boundary rigidity problem
for simple 2D manifolds,  injectivity of the solenoidal
part of the tensor field of order two was proven in \cite{Sh-2d}.
In \cite{SU-rig}, we considered $I$ on all geodesics and proved that the set of simple metrics on a fixed  manifold for which $I$  is s-injective is generic in $C^k(M)$, $k\gg1$. Previous results include s-injectivity for simple manifolds with curvature satisfying some explicit upper bounds \cite{Sh,Sh-sib,Pe}. A recent result by Dairbekov \cite{D} proves s-injectivity for non-trapping manifolds (not-necessarily convex) satisfying similar bounds, that in particular prevent the existence of  conjugate points.

Fix  another compact manifold $M_1$ with boundary such that $M_1^\text{\rm int}\supset M$, where $M_1^\text{\rm int}$ stands for the interior of $M_1$. Such a manifold is easy to construct in local charts, then glued together. 

\begin{definition}   \label{def_ms}
We say that the $C^k(M)$ (or analytic) metric $g$ on $M$ is \textbf{regular}, if $g$ has a $C^k$ (or analytic, respectively) extension  on $M_1$, such that for any $(x,\xi)\in T^*M$ there exists $\theta\in T_xM\setminus 0$ with $\langle \xi, \theta\rangle =0$ such that there is a geodesic segment $\gamma_{x,\theta}$ through  $(x,\theta)$ such that

(a) the endpoints of $\gamma_{x,\theta}$ are in  $\Mint_1\setminus M$.

(b)  there are no conjugate points on $\gamma_{x,\theta}$.

\noindent
Any geodesic satisfying (a), (b) is called a \textbf{simple}  geodesic. 
\end{definition}

Note that we allow the geodesics in $\Gamma$ to self-intersect. 

Since we do not assume that $M$ is convex, given $(x,\theta)$ there might be two or more geodesic segments $\gamma_j$ issued from $(x,\theta)$ such that $\gamma_j\cap M$  have different numbers of connected components. Some of them might be simple, others might be not. For example for a kidney-shaped domain and a fixed $(x,\theta)$ we may have such segments so that the intersection with $M$ has only one, or two connected components. Depending on which point in $T^*M$ we target to recover the singularities, we may need the first, or the second extension. So simple geodesic segments through some $x$ (that we call simple geodesics through $x$) are uniquely determined by an initial point $x$ and a direction $\theta$ and its endpoints. In case of simple manifolds, the endpoints (of the only connected component in $M$, unless the geodesics does not intersect $M$) are not needed, they are a function of $(x,\theta)$. Another way to determine a simple geodesic is by parametrizing it with $(x,\eta)\in T(\Mint_1 \setminus M)$, such that $\exp_x{\eta}\in \Mint_1\setminus M$  then
\be{01}
\gamma_{x,\eta} = 
\left\{  \exp_x(t\eta), 0\le t\le1\right\}.
\ee
This parametrization induces a topology  on the set $\Gamma$ of simple geodesics through  points of $\Mint_1$.

\begin{definition}   \label{def_complete}
The set $\Gamma$ of geodesics is called \textbf{complete}, if 

(a)   $\forall (x,\xi)\in T^*M$ there exists a simple geodesic $\gamma\in \Gamma$ through $x$ such that $\dot\gamma$ is normal to $\xi$ at $x$.

(b) $\Gamma$ is open.
\end{definition}

In other words, a regular metric $g$ is a metric for which a complete set of  geodesics exists. 
Another way to express (a) is to say that
\be{N*}
N^*\Gamma := \left\{N^*\gamma;\; \gamma\in \Gamma\right\} \supset T^*M,
\ee
where $N^*\gamma$ stands for the conormal bundle of $\gamma$.

We always assume that all tensor fields defined in $M$ are extended as $0$ to $M_1\setminus M$. Notice that $If$ does not change if we replace $M$ by another manifold $M_{1/2}$ close enough to $M$ such that $M\subset M_{1/2}\subset M_1$ but keep $f$ supported in $M$. Therefore, assuming that $M$ has an analytic structure as before, we can always extend $M$ a bit to make the boundary analytic and this would keep $(M,\bo,g)$ regular. Then s-injectivity in the extended $M$ would imply the same in the original $M$, see \cite[Prop.~4.3]{SU-rig}. So from now on, we will assume that $(M,\bo)$ is analytic but $g$ does not need to be analytic. To define correctly a norm in $C^K(M)$, respectively $C^k(M_1)$, we fix a finite analytic atlas.

The motivation behind Definitions~\ref{def_ms},~\ref{def_complete}  is the following: if $g$ is regular, and  $\Gamma$ is any  complete set of  geodesics, we will show that  $I_\Gamma f=0$ implies that $f^s\in C^l(M)$,   where $l=l(k)\to\infty$, as $k\to\infty$, in other words, the so restricted X-ray transform resolves the singularities. 

The condition of $g$ being regular is an open one for $g\in C^k(M)$, i.e., it defines an open set.  Any simple metric on $M$ is regular but the class of regular metrics is substantially larger if $\dim M\ge3$ and allows manifolds not necessarily diffeomorphic to a ball. For regular metrics on $M$,  we do not impose convexity assumptions on the boundary; conjugate points are allowed as far as the metric is regular; $M$ does not need to be non-trapping. In two dimensions, a regular metric can not have conjugate points in $M$ but the class is still larger than that of simple metrics because we do not require strong convexity of $\partial M$. 

\medskip
\paragraph{\bf Example 1.} To construct a manifold with a regular metric $g$ that has conjugate points, let us start with a manifold of dimension at least three with at least one pair of conjugate points $u$ and $v$ on a geodesic $[a,b]\ni t\mapsto \gamma(t)$. We assume that $\gamma$ is non-selfintersecting.  Then we will construct $M$ as a tubular neighborhood of $\gamma$. For any $x_0\in\gamma$, define $S_{x_0} = \exp_{x_0}\{v;\; 
\langle v, \dot\gamma(x_0)\rangle=0,\; |v|\le\eps  \}$, and 
$M := \cup_{x_0\in\gamma} S_{x_0}$ with $\eps\ll1$. Then there are no conjugate points along the geodesics that can be loosely described as those ``almost perpendicular'' to $\gamma$ but not necessarily intersecting $\gamma$; and the union of their conormal bundles covers  $T^*M$. More precisely, fix $x\in M$, then $x\in S_{x_0}$ for some $x_0\in\gamma$. Let $0\not=\xi\in T^*_xM$. Then there exists $0\not=v\in T_xM$ that is both tangent to $S_{x_0}$ and normal to $\xi$. The geodesic through $(x,v)$ is then a simple one for $\eps\ll1$, and the latter can be chosen in a uniform way independent of $x$. To obtain a smooth boundary, one can perturb $M$ so that the new manifold is still regular. 

\medskip
\paragraph{\bf Example 2.} This is similar to the example above but we consider a neighborhood of a periodic trajectory. Let $M =\left\{(x^1)^2+(x^2)^2\le 1\right\}\times S^1$ be the interior of the torus in $\R^3$, with the flat metric $(dx^1)^2+(dx^2)^2+d\theta^2$, where $\theta$ is the natural coordinate on $S^1$ with period $2\pi$. All geodesics perpendicular to $\theta=\mbox{const.}$ are periodic. All geodesics perpendicular to them have lengths not exceeding $2$ and their conormal bundles cover the entire $T^*M$ (to cover the boundary points, we do need to extend the geodesics in a neighborhood of $M$). Then $M$ is a regular manifold that is trapping, and one can easily show that a small enough perturbation of $M$ is also regular, and may still be trapping. 

\medskip

The examples above are partial cases of a more general one. Let $(M',\partial M')$ be a simple compact Riemannian  manifold with boundary with $\dim M'\ge2$, and let $M''$ be a Riemannian compact manifold with or without boundary. Let $M$ be a small enough perturbation of $M'\times M''$. Then $M$ is regular.

Let $g$ be a fixed regular metric on $M$. The property of $\gamma$ being simple is stable under small perturbations. The parametrization by $(x,\eta)$ as in \r{01} clearly has two more dimensions that what is needed to determine uniquely $\gamma|_M$. Indeed, a parallel transport of $(x,\eta)$ along $\gamma_{x,\eta}$, close enough to $x$, will not change $\gamma|_M$, similarly, we can replace $\eta$ by $(1+\eps)\eta$, $|\eps|\ll1$.

We assume throughout this paper that $M$ satisfies the following.

\medskip
\textbf{Topological Condition:} Any path in $M$ connecting two boundary points  is homotopic to a polygon $c_1\cup \gamma_1 \cup c_2\cup\gamma_2\cup\dots \cup\gamma_k \cup c_{k+1}$  with the properties:

(i) $c_j$ are paths on $\bo$;

(ii) For any $j$, $\gamma_j =\tilde\gamma_j|_M$ for some $\tilde\gamma_j \in\Gamma$; $\gamma_j$ lie in $\Mint$ with the exception of its endpoints and is transversal to $\bo$ at both ends.

\medskip

\begin{theorem} \label{thm_an} \ 
Let $g$ be an analytic, regular metric on $M$. Let $\Gamma$ be a complete complex of  geodesics. Then $I_\Gamma$ is s-injective.
\end{theorem}

The proof is based on using analytic pseudo-differential calculus, see \cite{Sj-Ast, T}. This has been used before in integral geometry, see e.g., \cite{BQ, Q}, see also \cite{SU-rig}.

To formulate a stability estimate, we will parametrize the simple geodesics in a way that will remove the extra two  dimensions. Let $H_m$ be a finite collection of smooth  hypersurfaces in $\Mint_1$. Let $\mathcal{H}_m$ be an open subset of $\{(z,\theta)\in SM_1; \; z\in H_m, \theta\not\in T_zH_m \}$,  and let $\pm l_m^\pm(z,\theta)\ge0$ be two continuous  functions. Let $\Gamma(\mathcal{H}_m)$ be the set of geodesics 
\be{5}
\Gamma(\mathcal{H}_m) = \left\{\gamma_{z,\theta}(t); \; l_m^-(z,\theta)\le t\le l_m^+(z,\theta), \; (z,\theta)\in \mathcal{H}_m  \right\},
\ee
that, depending on the context, is considered either as a family of curves, or as a point set. We also assume that each $\gamma\in \Gamma(\mathcal{H}_m)$ is a simple geodesic.

If $g$ is simple, then one can take a single $H=\partial M_1$ with $l^-=0$ and an appropriate $l^+(z,\theta)$. If $g$ is regular only, and $\Gamma$ is any  complete set of geodesics, then any small enough neighborhood of a simple geodesic in $\Gamma$ has the properties listed above and by a compactness argument on can choose a finite complete set of such $\Gamma(\mathcal{H}_m)$'s, that is included in the original $\Gamma$, see Lemma~\ref{lemma_H}.

Given $\mathcal{H}=\{\mathcal{H}_m\}$ as above, we consider an open set $\mathcal{H'}=\{\mathcal{H}_m'\}$, such that $\mathcal{H}_m' \Subset \mathcal{H}_m$, and let $\Gamma(\mathcal{H}_m')$ be the associated set of geodesics defined as in \r{5}, with the same $l_m^\pm$. Set $\Gamma(\mathcal{H})=\cup\Gamma(\mathcal{H}_m)$, $\Gamma(\mathcal{H}')=\cup\Gamma(\mathcal{H}_m')$. 

The restriction $\gamma\in \Gamma(\mathcal{H}_m')\subset \Gamma(\mathcal{H}_m)$ can be modeled by introducing a weight function $\alpha_m$ in $\mathcal{H}_m$, such that $\alpha_m=1$ on $\mathcal{H}_m'$, and $\alpha_m=0$ otherwise. More generally, we allow $\alpha_m$ to be smooth but still supported in $\mathcal{H}_m$. We then write $\alpha=\{\alpha_m\}$, and we say that $\alpha\in C^k(\mathcal{H})$, if $\alpha_m\in C^k(\mathcal{H}_m)$, $\forall m$. 

We consider $I_{\alpha_m}=\alpha_mI$, or more precisely, in the coordinates $(z,\theta) \in \mathcal{H}_m$,
\be{I_a0}
I_{\alpha_m}f = \alpha_m(z,\theta) \int_0^{l_m(z,\theta)} \big \langle  f(\gamma_{z,\theta}), \dot \gamma_{z,\theta}^2\big \rangle \,\d t, \quad (z,\theta)\in  \mathcal{H}_m.
\ee
Next, we set
\be{Na}
I_\alpha = \{ I_{\alpha_m} \}, \quad N_{\alpha_m} = I_{\alpha_m}^* I_{\alpha_m} = I^*|\alpha_m|^2I, \quad 
N_\alpha = \sum N_{\alpha_m},
\ee
where the adjoint is taken w.r.t.\ the measure $\d\mu := |\langle \nu(z),\theta\rangle | \,\d S_z\,\d \theta$ on $\mathcal{H}_m$,  $\d S_z\,\d \theta$ being the induced measure on $SM$, and $\nu(z)$ being a unit normal to $H_m$.

S-injectivity of  $N_\alpha$ is equivalent to s-injectivity for $I_\alpha$, which in turn is equivalent to s-injectivity of $I$ restricted to $\supp\alpha$, see  Lemma~\ref{lemma_1}. The space $\tilde{H}^2$ is defined in Section~\ref{sec_prel}, see \r{S24}.

\begin{theorem}   \label{thm_stab} \ 

(a) Let $g=g_0\in C^k$, $k\gg1$ be regular, and let $\mathcal{H}'\Subset\mathcal{H}$ be as above with $\Gamma(\mathcal{H}')$ complete.   Fix $\alpha = \{\alpha_m\}\in C^\infty$ with  $\mathcal{H}_m'  \subset\supp\alpha_m\subset \mathcal{H}_m$. 
Then if $I_\alpha$ is s-injective, we have 
\be{est}
\|f^s\|_{L^2(M)}  \le C
\|N_{\alpha} f\|_{\tilde H^2(M_1)}.
\ee

(b) Assume that $\alpha=\alpha_g$ in (a)  depends on $g \in C^k$, so that $C^k(M_1) \ni g \to C^l(\mathcal{H}) \ni \alpha_g$ is continuous with  $l\gg1$,  $k\gg1$.  Assume that $I_{g_0,\alpha_{g_0}}$ is s-injective. Then estimate \r{est} remains true for $g$ in a small enough neighborhood of $g_0$ in $C^k(M_1)$ with a uniform constant $C>0$. 
\end{theorem}

In particular, Theorem~\ref{thm_stab} proves a locally uniform stability estimate for the class of non-trapping manifolds considered in \cite{D}.

Theorems~\ref{thm_an}, \ref{thm_stab} allow us to formulate  generic uniqueness results. One of them is formulated below. Given a family of metrics $\mathcal{G}\subset C^k(M_1)$, and $U_g\subset T(\Mint_1\setminus M)$, depending on the metric $g\in \mathcal{G}$, we say that $U_g$ depends continuously on $g$, if for any $g_0\in \mathcal{G}$, and any compact $K\subset U^\text{int}_{g_0}$, we have $K\subset U^\text{int}_{g}$ for $g$ in a small enough neighborhood of $g_0$ in $C^k$. In the next theorem, we take $U_g=\Gamma_g$, that is identified with the corresponding set of $(x,\eta)$ as in \eqref{01}.

\begin{theorem}  \label{thm_I} 
Let $\mathcal{G}\subset C^k(M_1)$ be an open set of regular metrics on $M$, and let for  each $g\in\mathcal{G}$, $\Gamma_g$ be a   complete set of  geodesics related to $g$ and continuously depending on $g$. Then for $k\gg0$, there is an open and dense subset  $\mathcal{G}_s$ of $\mathcal{G}$, such that the corresponding X-ray transform $I_{\Gamma_g}$ is s-injective.
\end{theorem}

Of course, the set $\mathcal{G}_s$ includes all real analytic metrics in $\mathcal{G}$. 

\begin{corollary}  \label{cor_1} 
Let $\mathcal{R}(M)$ be the set of all regular $C^k$ metrics on $M$ equipped with the $C^k(M_1)$ topology. Then for $k\gg1$, the subset of metrics for which the   X-ray transform $I$ over all simple geodesics  is s-injective, is open and dense in  $\mathcal{R}(M)$.
\end{corollary}

The results above extend the generic results in \cite{SU-rig}, see also \cite{SU-Duke}, in several directions: the topology of $M$ may not be trivial, we allow conjugate points but we use only geodesics without conjugate  points; the boundary does not need to be convex; and we use incomplete data, i.e., we use integrals over subsets of geodesics only.

In Section~\ref{sec_f}, we  discuss  versions of those results for the X-ray transform of vector fields and functions, where the proofs can be simplified. Our results remain true for tensors of any order $m$, the necessary modifications are addressed in the key points of our exposition. To keep the paper readable, we restrict ourselves to orders $m=2,1,0$.

\section{Preliminaries}  \label{sec_prel}

We say that $f$ is analytic in some subset $U$ of an analytic manifold, not necessarily open, if $f$ can be extended analytically to some open set containing $U$. Then we write $f\in \mathcal{A}(U)$. 
Let $g\in \mathrm{C^k}(M)$, $k\gg2$ or $g\in\mathcal{A}(M)$ be a Riemannian metric in $M$. We work with symmetric 2-tensors $f=\{f_{ij}\}$ and with 1-tenors/differential forms $v_j$ (the notation here and below is in any local coordinates). We use freely the Einstein summation convention and the convention for raising and lowering indices. We think of $f_{ij}$ and  $f^{ij}= f_{kl}g^{ki}g^{lj}$ as different representations of the same tensor. If $\xi$ is a covector at $x$, then its components are denoted by   $\xi_j$, while $\xi^j$ is defined as $\xi^i = g^{ij}\xi_j$. Next, we denote $|\xi|^2=\xi_i\xi^i$, similarly for vectors that we usually denote by $\theta$.  If $\theta_1$, $\theta_2$ are two vectors, then $\langle \theta_1,\theta_2\rangle$ is their inner product. If $\xi$ is a covector, and $\theta$ is a vector, then $\langle \xi,\theta\rangle$ stands for $\xi(\theta)$. This  notation choice is partly justified by identifying $\xi$ with a vector, as above.

The geodesics of $g$ can be also viewed as the $x$-projections of the bicharacteristics of the Hamiltonian $E_g(x,\xi)=\frac12  g^{ij}(x) \xi_i\xi_j$. The energy level $E_g=1/2$ corresponds to parametrization with the arc-length parameter. For any geodesic $\gamma$, we have $f^{ij}(x)\xi_i\xi_j= f_{ij}(\gamma(x)) \dot\gamma^i(t) \dot\gamma^j(t)$, where $(x,\xi) = (x(t),\xi(t))$ is the bicharacteristic with $x$-projection equal to $\gamma$.

\subsection{Semigeodesic coordinates near a simple geodesic and boundary normal coordinates.}  \label{sec_sm} 

Let $[l^-,l^+]\allowbreak \ni t\mapsto  \gamma_{x_0,\theta_0}(t)$ be a simple geodesic through $x_0 =\gamma_{x_0,\theta_0}(0)  \in M_1$ with $\theta_0\in S_{x_0}M_1$. 
The map $t\theta\mapsto \exp_{x_0}(t\theta)$ is a local diffeomorphism for $\theta$ close enough to $\theta_0$ and $t\in[l^-,l^+]$ by our simplicity assumption but may not be a global one, since $\gamma_{x_0,\theta_0}$ may self-intersect. On the other hand, there can be finitely many intersections only and we can assume that each subsequent intersection happens on a different copy of $M$. In other words, we think of $\gamma_0$ as belonging to a new manifold that is a small enough neighborhood of $\gamma_0$, and there are no self-intersections there. The local charts of that manifold are defined through the  exponential map above. Therefore,  when working near $\gamma_{x_0,\theta_0}$ we can assume that $\gamma_{x_0,\theta_0}$ does not intersect itself. We will use this in the proof of Proposition~\ref{lemma_wf}. Then one can choose a  neighborhood $U$ of  $\gamma_0$ and  normal   coordinates centered at $x_0$ there, denoted by $x$ again, such that  
the radial lines $t\mapsto t\theta$, $\theta=\text{const.}$, are geodesics. If $g\in C^k$, then we lose two derivatives and the new metric is in $C^{k-2}$; if $g$ is analytic near $\gamma_0$, then the coordinate change can be chosen to be analytic, as well. 

If in the situation above, let $x_0\not\in M$, and moreover, assume that the part of $\gamma_{x_0,\theta_0}$ corresponding to $t<0$ is still outside $M$. Then, one can consider $(\theta,t)$ as polar coordinates on $T_{x_0}M$. Considering them as Cartesian coordinates there, see also \cite[sec.~9]{SU-Duke}, one gets coordinates  $(x',x^n)$ near $\gamma_{x_0,\theta_0}$ so that the latter is given by $\{(0,\dots,0,t), \; 0\le t\le l^+\}$, $g_{in}=\delta_{in}$, and $\Gamma_{nn}^i=\Gamma_{in}^n=0$, $\forall i$. Given $x\in \R^n$, we write $x' = (x^1,\dots,x^{n-1})$. Moreover, the lines $x'=\text{const.}$, $|x'|\ll1$, $x^n=t\in[0,l^+]$ are geodesics in $\Gamma$, as well. We will call those coordinates  semigeodesic coordinates near $\gamma_{x_0,\theta_0}$.

We will often use boundary normal (semi-geodesic) coordinates  $(x',x^n)$ near a boundary point. If $x'\in\R^{n-1}$ are local coordinates on $\bo$, and $\nu(x')$ is the interior unit normal, for $p\in M$ close enough to $\bo$, they are defined by $\exp_{(x',0)}x^n\nu=p$. Then  $x^n=0$ defines $\partial M$, $x^n>0$ in $M$, $x^n = \text{dist}(x,\partial M)$. The metric $g$ in those coordinates again satisfies  $g_{in}=\delta_{in}$, and $\Gamma_{nn}^i=\Gamma_{in}^n=0$, $\forall i$. We also use the  convention that  all Greek indices take values from $1$ to $n-1$. 
In fact, the semigeodesic coordinates in the previous paragraph  are boundary normal coordinates to a  small part of the geodesic ball centered at $x_0= \gamma_{x_0,\theta_0}(0)$ with radius $\eps$, $0<\eps\ll1$.

\subsection{Integral representation of the normal operator.}   \label{sec_int}

We define the $L^2$ space of symmetric tensors $f$ with inner product
\[
(f,h) =  \int_M  \langle f, \bar h\rangle (\det g)^{1/2}\,\d x,
\]
where, in local coordinates, $\langle f, \bar h\rangle = f_{ij}\bar h^{ij}$. Similarly, we define the $L^2$ space of 1-tensors (vector fields, that we identify with 1-forms) and the $L^2$ space of functions in $M$. Also, we will work in Sobolev $H^s$ spaces of 2-tensors, 1-forms and functions. In order to keep the notation simple, we will use the same notation $L^2$ (or $H^s$) for all those  spaces and it will be clear from the context which one we mean.

In the fixed finite atlas on $M$, extended to $M_1$, the norms $\|f\|_{C^k}$ and the $H^s$ norms below are correctly defined. In the proof, we will work in finitely many coordinate charts because of the compactness of $M$, and this justifies the equivalence of the correspondent $C^k$ and $H^s$ norms. 

We define the Hilbert space  $\tilde{H}^2(M_1)$ used in Theorem~\ref{thm_stab} as in \cite{SU-Duke,SU-rig}. Let $x=(x',x^n)$ be local coordinates in a neighborhood $U$ of a point on $\partial M$ such that $x^n=0$ defines $\partial M$. Then we  set
\[
\|f\|^2_{\tilde{H}^1(U)} = \int_U \Big(\sum_{j=1}^{n-1} |\partial_{x^j}f|^2 +|x^n\partial_{x^n}f|^2+|f|^2\Big)\, \d x.
\]
This can be extended to a small enough neighborhood $V$ of $\partial M$ contained in $M_1$. Then we set
\begin{equation}  \label{S24}
\|f\|_{\tilde{H}^2(M_1)} = \sum_{j=1}^{n}  \|\partial_{x^j}f\|_{\tilde{H}^1(V)} + \|f\|_{\tilde{H}^1(M_1)}.
\end{equation}
The space $\tilde{{H}}^2(M_1)$ has the property that for each $f\in H^1(M)$ (extended as zero outside $M$), we have $N f \in \tilde{H}^2(M_1)$. This is not true if we replace $\tilde{H}^2(M_1)$ by $H^2(M_1)$.

\begin{lemma}  \label{lemma_H}
Let $\Gamma_g$ and $\mathcal{G}$ be as in Theorem~\ref{thm_I}.  Then for $k\gg1$, for any $g_0\in\mathcal{G}$, there exist  $ \mathcal{H}'=\{\mathcal{H}_m'\}\Subset \mathcal{H}=\{\mathcal{H}_m\}$ such that $\Gamma(\mathcal{H})\Subset\Gamma_{g_0}$, and $\mathcal{H}'$, $\mathcal{H}$ satisfy the assumptions of Theorem~\ref{thm_stab}. Moreover, $\mathcal{H}'$ and $\mathcal{H}$ satisfy the assumptions of Theorem~\ref{thm_stab} for $g$ in a small enough neighborhood of $g_0$ in $C^k$. 
\end{lemma}

\begin{proof} 
Fix $g_0\in\mathcal{G}$ first. 
Given $(x_0,\xi_0)\in T^*M$, there is a simple geodesic $\gamma: [l^-, l^+] \to M_1$ in $\Gamma_{g_0}$  through $x_0$  normal to $\xi_0$ at $x_0$. Choose a small enough hypersurface $H$ through $x_0$ transversal to $\gamma\in \Gamma_{g_0}$, and local coordinates near $x_0$ as in Section~\ref{sec_sm} above,  so that $x_0=0$, $H$ is given by $x^n=0$, $\dot\gamma(0)=(0,\dots,0,1)$. Then one can set $\mathcal{H}_0 = \{x;\; x^n=0; \; |x'|<\eps\} \times\{ \theta;\; |\theta'|<\varepsilon  \}$, and $\mathcal{H}_0'$ is defined in the same way by replacing $\eps$ by $\eps/2$. We define $\Gamma(\mathcal{H}_0)$ as in \r{5} with $l^\pm(z,\theta)=l^\pm$. 
Then the properties required for $\mathcal{H}_0$, including the simplicity assumption are satisfied when $0<\eps\ll1$. Choose such an $\eps$, and replace it with a smaller one so that those properties are preserved under a small perturbation of $g$.   Any point in $SM$ close enough to $(x_0,\xi_0)$ still has a geodesic in $\Gamma(\mathcal{H}_0')$ normal to it. By a compactness argument, one can find a finite number of $\mathcal{H}_m'$ so that the corresponding $\Gamma(\mathcal{H}') = \cup \Gamma(\mathcal{H}'_m)$ is complete. 

The continuity property of $\Gamma_g$ w.r.t.\ $g$ guarantees that the construction above is stable under a small perturbation of $g$.
\end{proof}

Similarly to \cite{SU-Duke}, one can see that the map $I_{\alpha_m}  : L^2(M) \to L^2(\mathcal{H}_m,\,\d\mu)$ defined by \r{I_a0} is bounded, and therefore the {\em normal}\/ operator $N_{\alpha_m}$ defined in \r{Na} is a well defined bounded operator on $L^2(M)$. Applying the same argument to $M_1$, we see that $N_{\alpha_m} : M \to M_1$ is also bounded. 
By \cite{SU-Duke}, at least when $f$ is supported in the local chart near $x_0=0$ above, and $x$ is close enough to $x_0$, 
\be{N0}
\left[N_{\alpha_m}f\right]^{i'j'} (x) =  \int_0^\infty \int_{S_x\Omega} |\alpha_m^\sharp(x,\theta)|^2  \theta^{i'}  \theta^{j'} f_{ij}( \gamma_{x,\theta}(t) ) \dot \gamma_{x,\theta}^i(t) \dot\gamma_{x,\theta}^j(t)\, \d \theta\, \d t,
\ee
where $|\alpha_m^\sharp(x,\theta)|^2 = |\tilde\alpha_m(x,\theta)|^2 + |\tilde\alpha_m(x,-\theta)|^2$, and $\tilde\alpha_m$ is the extension of $\alpha_m$ as constant along the geodesic through $(x,\theta)\in \mathcal{H}_m$; and equal to $0$ for all other points not covered by such geodesics. Formula \r{N0} has an invariant meaning and holds without the restriction on $\supp f$. On the other hand, if $\supp f$ is small enough (but not necessarily near $x_0$), $y=\exp_x(t\theta)$ defines a local diffeomorphism $t\theta\mapsto y\in \supp f$, therefore after making the change of variables $y=\exp_x(t\theta)$, see \cite{SU-Duke}, this becomes
\be{N}
N_{\alpha_m}f (x) = 
\frac1{\sqrt{\det g}}
 \int A_m(x,y) \frac{f^{ij}(y)}{\rho(x,y)^{n-1}}
\frac{\partial\rho}{\partial y^i} \frac{\partial\rho}{\partial y^j}
\frac{\partial\rho}{\partial x^k} \frac{\partial\rho}{\partial x^l}\,
\!\det\frac{\partial^2(\rho^2/2)}{\partial x\partial y} \,\d y, 
\ee
where
\be{A}
A_m(x,y) = \big|\alpha_m^\sharp\!\left(x, \text{grad}_x \rho(x,y)\right)\big|^2,
\ee
$y$ are any local coordinates near $\supp f$, and $\rho(x,y) = |\exp^{-1}_x y|$. Formula \r{N}  can be also understood invariantly by considering $\d_x \rho$ and $\d_y\rho$ as tensors. 
For arbitrary $f\in L^2(M)$ we  use a partition of unity in $T\Mint_1$ to express $N_{\alpha_m}f(x)$ as a finite sum of integrals as above, for $x$ near any fixed $x_0$.

We get in particular that  $N_{\alpha_m}$ has the pseudolocal  property, i.e., its Schwartz kernel is smooth outside the diagonal. As we will show below, similarly to the analysis in \cite{SU-Duke, SU-rig}, $N_{\alpha_m}$ is a \PDO\ of order $-1$. 

We always extend functions or tensors  defined in $M$ as $0$ outside $M$.  Then $N_\alpha f$ is well defined near $M$ as well and remains unchanged if $M$ is extended such that it is still in $M_1$, and $f$ is kept fixed.

\subsection{Decomposition of symmetric tensors.}  
For more details about the decomposition below, we refer to \cite{Sh}. 
Given a symmetric  2-tensor $f= f_{ij}$, we define the 1-tensor $\delta f$ called {\em divergence}  of $f$ by 
$$
[\delta f]_i = g^{jk} \nabla_k f_{ij},
$$ 
in any local coordinates, where $\nabla_k$ are the covariant derivatives of the tensor $f$. Given an 1-tensor (a vector field or an 1-form) $v$, we denote by $dv$ the 2-tensor called symmetric differential of $v$:
$$
[d v]_{ij} = \frac12\left(\nabla_iv_j+ \nabla_jv_i  \right).
$$
Operators $d$ and $-\delta$ are formally adjoint to each other 
in $L^2(M)$.  
It is easy to see that for each smooth $v$ with $v=0$ on $\partial M$, we have $I(d v)(\gamma)=0$ for any geodesic $\gamma$ with endpoints on $\bo$. This follows from the identity  
\be{v}
\frac{\d }{\d  t}  \langle v(\gamma(t)), \dot\gamma(t) 
\rangle \allowbreak  = \allowbreak \langle
dv(\gamma(t)),  \dot\gamma^2(t) \rangle.
\ee 
If $\alpha=\{\alpha_m\}$ is as in the Introduction, we get
\be{dv}
I_\alpha  (dv)=0, \quad \forall v\in C_0^1(M),
\ee
and this can be extended to $v\in H_0^1(M)$ by continuity. 

It is known (see \cite{Sh} and  \r{10} below) that for $g$ smooth enough, each symmetric tensor $f\in L^2(M)$ admits unique orthogonal decomposition $f=f^s+d v$ into a {\em solenoidal}\/ tensor  $\mathcal{S}f :=f^s $ and a {\em potential}\/ tensor $\mathcal{P}f :=d v$, such that both terms are in $L^2(M)$, $f^s$ is solenoidal, i.e., $\delta f^s=0$ in $M$, and $v\in H^1_0(M)$ (i.e., $v=0$ on $\partial M$). In order to construct this decomposition, introduce the operator $\upDelta^s = \delta d$ acting on vector fields. This operator is elliptic in $M$, and the Dirichlet problem satisfies the Lopatinskii condition. Denote by $\upDelta^s_D$ the Dirichlet realization of $\upDelta^s$ in $M$. Then
\begin{equation}  \label{9}
v = \left(\upDelta^s_D\right)^{-1}\delta f, \quad
f^s = f - d \left(\upDelta^s_D\right)^{-1}\delta f.
\end{equation} 
Therefore, we have
$$
\mathcal{P} = d \left(\upDelta^s_D\right)^{-1}\delta, \quad \mathcal{S} = \Id-\mathcal{P},
$$
and for any $g \in C^1(M)$, the maps
\be{10}
(\upDelta^s_D)^{-1}: H^{-1}(M) \to H_0^{1}(M),  \quad 
\mathcal{P},  \mathcal{S} :  L^2(M) \longrightarrow L^2(M)
\ee
are bounded and depend continuously on $g$, see \cite[Lemma~1]{SU-rig} that easily generalizes for manifolds. 
This  admits the following easy generalization:  for $s=0,1,\dots$, the resolvent above also continuously maps $H^{s-1}$ into $H^{s+1} \cap H_0^1$, similarly, $\mathcal{P}$ and $\mathcal{S}$ are bounded in $H^{s}$, if $g\in C^k$, $k\gg1$ (depending on $s$).  Moreover those operators depend continuously on $g$. 

Notice that even when $f$ is smooth and $f=0$ on $\partial M$, then $f^s$ does not need to vanish on $\partial M$. In particular, $f^s$, extended as $0$ to $M_1$, may not be solenoidal anymore. 
To stress on the dependence on the manifold, when needed, we will use the notation $v_M$ and $f^s_M$ as well. 

Operators $\mathcal{S}$ and $\mathcal{P}$ are orthogonal projectors. 
The problem about the s-injectivity of $I_\alpha$ then can be posed as follows: if $I_\alpha f=0$, show that $f^s=0$, in other words, show that $I_\alpha$ is injective on the subspace $\mathcal{S}L^2$ of  solenoidal tensors. Note that by \r{dv} and \r{Na},
\be{11}
N_\alpha = N_\alpha\mathcal{S}=\mathcal{S} N_\alpha, \quad \mathcal{P} N_\alpha=N_\alpha\mathcal{P}=0.
\ee

\begin{lemma}  \label{lemma_1} Let $\alpha=\{\alpha_m\}$ with $\alpha_m\in C_0^\infty(\mathcal{H}_m)$ be as in the Introduction. The following statements are equivalent:

(a) $I_\alpha$ is s-injective on $L^2(M)$;

(b) $N_\alpha : L^2(M) \to L^2(M)$ is s-injective;

(c) $N_\alpha : L^2(M) \to L^2(M_1)$ is s-injective;

(d) If $\Gamma^\alpha_m$ is the set of geodesics issued from $(\supp\alpha_m)^\text{\rm int}$ as in \r{5}, and  $\Gamma^\alpha = \cup\Gamma_m^\alpha$, then $I_{\Gamma^\alpha}$ is s-injective.
\end{lemma}

\begin{proof} 
Let $I_\alpha$ be s-injective, and assume that $N_\alpha  f=0$ in $M$ for some $f\in L^2(M)$. Then 
$$
0 = (N_\alpha   f,f)_{L^2(M)} = \sum \|\alpha_m I f\|_{L^2(\mathcal{H}_m,\d\mu)}^2 \quad \Longrightarrow \quad f^s=0.
$$
This proves the implication $(a) \Rightarrow (b)$. Next, $(b) \Rightarrow (c)$ is immediate.  
Assume (c) and let $f\in L^2(M)$ be such that $I_\alpha f=0$. Then $N_\alpha f=0$ in $M_1$, therefore $f^s=0$. Therefore,  $(c) \Rightarrow (a)$. Finally,  $(a) \Leftrightarrow (d)$ follows directly form the definition of $I_\alpha$. 
\end{proof}

\paragraph{\bf Remark.} Lemma~\ref{lemma_1} above, and Lemma~\ref{lemma_bd}(a) in next section show that $(\supp\alpha_m)^\text{int}$ in (d) can be replaced by $\supp\alpha_m$ if  $\Gamma^\alpha$ is a complete  set of  geodesics.

\section{Microlocal Parametrix of $N_\alpha$}

\begin{proposition}  \label{pr_2}\ 
Let $g=g_0\in C^k(M)$ be a regular metric on $M$, and let $\mathcal{H}'\Subset \mathcal{H}$ be as in Theorem~\ref{thm_stab}. 

(a)   Let $\alpha$ be as in Theorem~\ref{thm_stab}(a). 
Then for any $t=1,2,\dots$, there exists $k>0$ and a bounded linear operator
$$
Q : \tilde H^2(M_1)  \longmapsto \mathcal{S} L^2(M),
$$
such that
\be{F}
QN_\alpha f = f_M^s +Kf, \quad \forall f\in H^1(M),
\ee
where $K :H^1(M) \to \mathcal{S} H^{1+t}(M)$ extends to $K :L^2(M) \to \mathcal{S} H^{t}(M)$. If $t=\infty$, then $k=\infty$. 

(b) Let $\alpha =\alpha_g$ be as in Theorem~\ref{thm_stab}(b). 
Then, for $g$ in some $C^k$ neighborhood of $g_0$,  (a) still  holds  and $Q$ can be constructed so that $K$ would depend continuously on $g$.
\end{proposition}

\begin{proof} 
A brief sketch of our proof is the following: We construct first  a parametrix that  recovers microlocally $f^s_{M_1}$ from $N_\alpha f$. Next we will compose this parametrix with  the operator $f_{M_1}^s \mapsto f_M^s$ as in \cite{SU-Duke, SU-rig}. Part (b) is based on a perturbation argument for the Fredholm equation \r{F}. The need for such two step construction is due to the fact that in the definition of $f^s$, a solution to a certain boundary value problem is involved, therefore near $\bo$, our construction is not just a parametrix of a certain elliptic \PDO. This is the reason for losing one derivative in \r{est}. For tensors of orders 0 and 1, there is no such loss, see \cite{SU-Duke} and \r{est1}, \r{est2}.

As in \cite{SU-rig}, we will  work with \PDO s with symbols of finite smoothness $k\gg1$. All operations we are going to perform would require finitely many derivatives of the amplitude and finitely many seminorm estimates. In turn, this would be achieved if $g\in C^k$, $k\gg1$ and the corresponding \PDO s will depends continuously on $g$. 

Recall \cite{SU-Duke,SU-rig} that for simple metrics, $N$ is a \PDO\ in $\Mint$ of order $-1$ with principal symbol that is not elliptic but $N+|D|^{-1}\mathcal{P}$ is elliptic. This is a consequence of the following. We will say that $N_\alpha$ (and any other  \PDO\ acting on symmetric tensors) is {\em elliptic on solenoidal tensors}, if for any $(x,\xi)$, $\xi\not=0$,  $\sigma_p(N_\alpha)^{ijkl}(x,\xi)f_{kl}=0$ and $\xi^i f_{ij}=0$ imply $f=0$. Then $N$ is elliptic on solenoidal tensors, as shown in \cite{SU-Duke}. That definition is motivated by the fact that the principal symbol of $\delta$ is given by $f_{ij} \mapsto \mathrm{i}\xi^if_{ij}$, and s-injectivity is equivalent to the statement that $Nf=0$ and $\delta f=0$ in $M$ imply $f=0$. Note also that the principal symbol of $d$ is given by $v_j \mapsto (\xi_i v_j +\xi_j v_i)/2$, and $\sigma_p(N)$ vanishes on tensors represented by the r.h.s.\ of the latter. We will establish similar properties of $N_\alpha$ below.

Let $N_{\alpha_m}$ be as in Section~\ref{sec_int} with $m$ fixed.

\begin{lemma}  \label{lemma_2}
$N_{\alpha_m}$ is a classical \PDO\ of order $-1$ in $\Mint_1$. It is elliptic on solenoidal tensors at $(x_0,\xi^0)$ if and only if there exists $\theta_0\in T_{x_0}M_1\setminus 0$ with $\langle  \xi^0, \theta_0 \rangle=0$ such that  $\alpha_0(x_0,\theta_0)\not =0$.  
The principal symbol $\sigma_p(N_{\alpha_m})$ vanishes on tensors of the kind $f_{ij} = (\xi_i v_j +\xi_j v_i)/2$ and is non-negative on tensors satisfying $\xi^if_{ij}=0$. 
\end{lemma}

\begin{proof} We established the pseudolocal property already, and formulas \r{N0}, \r{N} together with the partition of unity argument following them imply that it is enough to work with $x$ in a small neighborhood of a fixed $x_0\in \Mint_1$, and with $f$ supported there as well. Then we work in local coordinates near $x_0$. To express $N_{\alpha_m}$ as a pseudo-differential operator, we proceed as in \cite{SU-Duke, SU-rig}, with a starting point \r{N}.  Recall that for $x$ close to $y$ we have
\[
\begin{split}
\rho^2(x,y)&=G^{(1)}_{ij}(x,y)(x-y)^i(x-y)^j,\\
\frac{\partial\rho^2(x,y)}{\partial x^j} &=2 G^{(2)}_{ij}(x,y)(x-y)^i,\\
\frac{\partial^2\rho^2(x,y)}{\partial x^j\partial y^j} &=2   G^{(3)}_{ij}(x,y),
\end{split}
\]
where $G^{(1)}_{ij}$, $G^{(2)}_{ij}$ $G^{(3)}_{ij}$   are smooth and on the diagonal. We have 
$$
G^{(1)}_{ij}(x,x)=  G^{(2)}_{ij}(x,x)= G^{(3)}_{ij}(x,x)=  g_{ij}(x).
$$
Then $N_{\alpha_m}$ is a pseudo-differential operator with amplitude
\begin{equation} \label{a20'}
\begin{split}  
M_{ijkl}(x,y,\xi) &
= \int e^{-\mathrm{i}\xi\cdot z}\left(G^{(1)}z\cdot z\right)^{\frac{-n+1}2-2}     \big|\alpha_m^\sharp(x,g^{-1}G^{(2)} z)\big|^2 \\       
&\qquad \times \big[G^{(2)}z\big]_i \big[G^{(2)}z\big]_j\big[\wtilde G^{(2)}z\big]_k
\big[\wtilde G^{(2)} z\big]_l
\frac{\det G^{(3)}}{\sqrt{\det g}}  \,\d z,
\end{split}
\end{equation}
where $\wtilde G^{(2)}_{ij}(x,y)=  G^{(2)}_{ij}(y,x)$. 
As in \cite{SU-rig}, we note that $M_{ijkl}$ is the Fourier transform of a positively homogeneous distribution in the $z$ variable, of order $n-1$. Therefore, $M_{ijkl}$ itself is positively homogeneous of order $-1$ in $\xi$. Write 
\be{m}
M(x,y,\xi) 
= \int e^{-\mathrm{i}\xi\cdot  z}|z|^{ -n+1} m(x,y,\theta) \,\d z, \quad \theta=z/|z|,
\ee
where
\be{30}
\begin{split}
m_{ijkl}(x,y,\theta) = &
\left(G^{(1)}\theta\cdot \theta\right)^{\frac{-n+1}2-2}  \big|\alpha_m^\sharp(x,g^{-1}G^{(2)} \theta)\big|^2 \\ & \times
\big[G^{(2)}\theta \big]_i \big[G^{(2)}\theta\big]_j\big[\wtilde G^{(2)}\theta\big]_k
\big[\wtilde G^{(2)} \theta\big]_l
\frac{\det G^{(3)}}{\sqrt{\det g(x)}},
\end{split}
\ee
and pass to polar coordinates $z=r\theta$. Since $m$ is an even function of $\theta$, smooth w.r.t.\ all variables,  we get (see also \cite[Theorem~7.1.24]{H})
\be{M}
M(x,y,\xi) = \pi \int_{|\theta|=1}  m(x,y,\theta)\delta(\theta\cdot\xi) \,\d \theta.
\ee
This proves that $M$ is an amplitude of order $-1$. 

To obtain the principal symbol, we set $x=y$ above (see also \cite[sec.~5]{SU-Duke} to get
\be{32}
\sigma_p(N_{\alpha_m}) (x,\xi) = M(x,x,\xi) = \pi \int_{|\theta|=1}  m(x,x,\theta)\delta(\theta\cdot\xi) \,\d \theta,
\ee
where
\be{33}
m^{ijkl}(x,x,\theta) =  \big|\alpha_m^\sharp(x,\theta)\big|^2 \sqrt{\det g(x)}\left(g_{ij}(x)\theta ^i \theta^j\right)^{\frac{-n+1}2-2}  \theta^i \theta^j  \theta^k  \theta^l  .
\ee

To prove ellipticity of $M(x,\xi)$ on solenoidal tensors  at $(x_0,\xi^0)$, notice that for any  symmetric real $f_{ij}$, we have
\be{33a}
m^{ijkl}(x_0,x_0,\theta) f_{ij} f_{kl} =  \big|\alpha_m^\sharp(x_0,\theta) \big|^2 \sqrt{\det g(x_0)}\left(g_{ij}(x_0)\theta ^i \theta^j\right)^{\frac{-n+1}2-2}  \!\left( f_{ij} \theta^i \theta^j \right)^2\ge0.
\ee
This, \r{32}, and the assumption $\alpha_m(x_0,\theta_0)\not= 0$ imply that $M^{ijkl}(x_0,x_0,\xi^0)f_{ij} f_{kl}=0$ yields $f_{ij}\theta^i \theta^j=0$ for  $\theta$ perpendicular to $\xi^0$, and close enough to $\theta_0$. If in addition $(\xi^0)^j f_{ij}=0$, then this implies $f_{ij}\theta^i \theta^j=0$ for $\theta\in \n(\theta_0)$, and that easily implies that it vanishes for all $\theta$. Since $f$ is symmetric, this means that $f=0$. 

The last statement of the lemma follows directly from \r{32}, \r{33}, \r{33a}.

Finally, we note that \r{33}, \r{33a} and  the proof above generalizes easily for tensors of any order. 
\end{proof}

We continue with the proof of Proposition~\ref{pr_2}. Since (b) implies (a), we will prove (b) directly. Notice that $\mathcal{H}'$ and $\mathcal{H}$ satisfy the properties listed in the Introduction, right before Theorem~\ref{thm_stab}, if $g=g_0$. On the other hand, those properties are stable under small $C^k$ perturbation of $g_0$. We will work here with metrics $g$ close enough to $g_0$.

By  Lemma~\ref{lemma_2}, since $\Gamma(\mathcal{H}')$ is complete, $N_\alpha$ defined by \r{Na} is elliptic on solenoidal tensors in $M$. 
The rest of the proof is identical to that of \cite[Proposition~4]{SU-rig}. We will give a  brief sketch of it. 
To use the ellipticity of $N_\alpha$  on solenoidal tensors, we complete $N_\alpha$ to an elliptic \PDO\ as in \cite{SU-rig}. Set
\be{W}
W = N_\alpha  +
|D|^{-1}\mathcal{P}_{M_1},
\ee
where $|D|^{-1}$ is a properly supported parametrix of $(-\Delta_g)^{1/2}$ in $\n(M_1)$. The resolvent $(-\Delta^s_{M_1,D})^{-1}$ involved in $\mathcal{P}_{M_1}$ and $\mathcal{S}_{M_1}$ can be expressed as $R_1+R_2$, where $R_1$ is any parametrix near $M_1$, and $R_2 : L^2_{\text{comp}}(M_1) \to C^l(M_1)$, $R_2: H^l(M_1) \to H^{l+2}(M_1)$, where $l=l(k)\gg1$, if $k\gg1$. Then $W$ is an elliptic \PDO\ inside $M_1$ of order $-1$ by Lemma~\ref{lemma_2}. 

Let $P$ be a properly supported parametrix for $W$  of finite order, i.e., $P$ is a classical \PDO\ in the interior of $M_1$ of order $1$ with amplitude of finite smoothness, such that 
\be{34}
PW=\Id +K_1,
\ee
and $K_1 : L^2_{\text{comp}} (M_1)\to H^l(M_1)$ with $l$ as above.  Then
$$
P_1 := \mathcal{S}_{M_1}P
$$
satisfies 
\be{35}
P_1 N_\alpha  = \mathcal{S}_{M_1}+K_2,
\ee
where $K_2$ has the same property as $K_1$. To see this, it is enough to apply $\mathcal{S}_{M_1}$ to the left and right of \r{34} and to use \r{11}. 

Next step is to construct an operator that recovers $f^s_{M}$, given $f^s_{M_1}$, and to apply it to $P_1 N_\alpha -K_2$. In order to do this, it is enough first to construct a map $P_2$  such that if  $f^s_{M_1}$ and $v_{M_1}$ are the solenoidal part and the potential, respectively, corresponding to $f\in L^2(M)$ extended as zero to $M_1\setminus M$, then $P_2  : f^s_{M_1}\mapsto \left. v_{M_1}\right|_{\partial M}$. 
This is done as in \cite{SU-Duke} and \cite[Proposition~4]{SU-rig}. We also have
$$
P_2P_1 : \wtilde H^2(M_1) \to H^{1/2}(\partial M). 
$$
Then we showed in \cite[Proposition~4]{SU-rig} that one can set
$$
Q  = (\Id +dRP_2)P_1,
$$
where $R : h\mapsto u$ is the Poisson operator  for the Dirichlet problem  $\upDelta^s u=0$ in $M$, $u|_{\partial M} =h$.  

As explained above, we work with finite asymptotic expansions that require finite number of derivatives on the amplitudes of our \PDO s. On the other hand,  these amplitudes depend continuously on $g\in C^k$, $k\gg1$. As a result, all operators above depend continuously on $g\in C^k$, $k\gg1$. 
\end{proof}

The first part of next lemma generalizes similar results in \cite[Thm~2]{SU-Duke}, \cite{ Ch, SSU} to the present situation. The second part  shows that $I_\Gamma f=0$ implies that a certain $\tilde f$, with the same solenoidal projection, is flat at $\partial M$. This $\tilde f$ is defined by the property \r{l1_2} below.

\begin{lemma}  \label{lemma_bd}
Let $g\in C^k(M)$ be a regular metric, and let $\Gamma$ be a complete  set of  geodesics.  Then

(a) $\Ker I_\Gamma \cap \mathcal{S}L^2(M)$
 is finite dimensional and included in $C^l(M)$ with $l=l(k)\to\infty$, as $k\to\infty$. 

(b) If $I_\Gamma f=0$ with $f\in L^2(M)$, then there exists a vector field $v\in C^l(M)$, with $v|_{\partial M}=0$ and  $l$ as above, such that for $\tilde f := f-dv$ we have
\be{l1_1}
\partial^\alpha \tilde f|_{\partial M} =0, \quad | \alpha|\le l,
\ee
and in boundary normal coordinates near any point on $\partial M$ we have
\be{l1_2}
\tilde f_{ni}=0,\quad \forall i.
\ee
\end{lemma}

\begin{proof} 
Part (a) follows directly from Proposition~\ref{pr_2}. 

Without loss of generality, we may assume that $M_1$ is  defined as $M_1 = \{x,\; \mbox{dist}(x,M)\le\epsilon\}$, with $\epsilon>0$ small enough. 
By Proposition~\ref{pr_2}, applied to $M_1$,
\be{11_3}
f_{M_1}^s \in C^l(M_1),
\ee
where $l\gg1$, if $k\gg1$.  

Let $x=(x',x^n)$ be boundary normal coordinates in a neighborhood of some boundary point. We recall how to construct $v$ defined in $M$ so that \r{l1_2} holds, see \cite{SU1} for a similar argument  for the non-linear boundary rigidity problem, and \cite{E,Sh-sib,SU-Duke,SU-rig} for the present one. 
The condition $(f-dv)_{in}=0$ is equivalent to 
\begin{equation}  \label{a1}
\nabla_n v_i+\nabla_i v_n= 2f_{in}, \quad v|_{x^n=0}=0, \quad i=1,\dots,n.
\end{equation}
Recall that $\nabla_i v_j = \partial_i v_j-\Gamma_{ij}^kv_k$, and that in those coordinates, $\Gamma_{nn}^k=\Gamma_{kn}^n=0$. If $i=n$, then \r{a1} reduces to  $\nabla_n v_n=\partial_n v_n=f_{nn}$, $v_n=0$ for $x^n=0$; we solve this by integration over $0\le x^n\le \eps \ll1$; this gives us $v_n$. Next, we solve the remaining linear system of $n-1$ equations for $i=1,\dots,n-1$ that is of the form $\nabla_nv_i=2f_{in}-\nabla_iv_n$, or, equivalently, 
\begin{equation}  \label{a1'}
\partial_n v_i-2\Gamma^\alpha_{ni}v_\alpha =  2f_{in}-\partial _iv_n, \quad v_i|_{x^n=0}=0, \quad i=1,\dots,n-1,
\end{equation}
(recall that $\alpha=1,\dots,n-1$). 
Clearly, if $g$ and $f$ are  smooth enough near $\partial M$, then so is $v$. If we set $f=f^s$ above (they both belong to $\Ker I_\Gamma$), then by (a) we get the statement about the smoothness of $v$. 
Since the condition \r{l1_2} has an  invariant meaning, this in fact defines a construction in some one-sided neighborhood of $\partial M$ in $M$. One can cut $v$ outside that neighborhood in a smooth way to define $v$ globally in $M$.  We also note that this can be done for tensors of any order $m$, see \cite{Sh-sib}, then we have to solve consecutively $m$  ODEs.

Let $\tilde f =f-dv$, where $v$ is as above. Then $\tilde f$ satisfies \r{l1_2}, and let
\be{37}
\tilde f^s_{M_1} = \tilde f - d\tilde v_{M_1}
\ee
be the solenoidal projection of $\tilde f$ in $M_1$. Recall that $\tilde f$,  according to our convention, is extended as zero in $M_1\setminus M$ that in principle, could create jumps across $\partial M$. Clearly, $\tilde f^s_{M_1} = f^s_{M_1}$ because $f-\tilde f=dv$ in $M$ with $v$ as in the previous paragraph, and this is also true in $M_1$ with $\tilde f$, $f$ and $v$ extended as zero (and then $v=0$ on $\bo_1$). In \r{37}, the l.h.s.\ is smooth in $M_1$ by \r{11_3}, and $\tilde f$ satisfies \r{l1_2} even outside $M$, where it is zero. Then one can get $\tilde v_{M_1}$ by solving \r{a1} with $M$ replaced by $M_1$, and $f$ there replaced by $\tilde f^s_{M_1}\in C^l(M_1)$. Therefore, one gets that $\tilde v_{M_1}$, and therefore $\tilde f$,  is smooth enough across $\partial M$, if $g\in C^k$, $k\gg1$, which proves \r{l1_1}. 

One can give the following alternative proof of \r{l1_1}: Let $N_\alpha$ be related to $\Gamma$, as in Theorem~\ref{thm_stab}. One can easily check that $N_\alpha$, restricted to tensors satisfying \r{l1_2}, is elliptic for $\xi_n\not=0$. Since $N_\alpha \tilde f=0$ near $M$, with $\tilde f$ extended as 0 outside $M$, as above, we get that this extension cannot have conormal singularities across $\bo$. This implies \r{l1_1}, at least when $g\in C^\infty$. The case of $g$ of finite smoothness can be treated by using parametrices of finite order in the conormal singularities calculus.
\end{proof}

\section{S-injectivity for analytic regular metrics} 
In this section, we prove Theorem~\ref{thm_an}. Let $g$ be an analytic regular metrics in $M$, and let $M_1\supset M$ be the manifold where $g$ is extended analytically according to Definition~\ref{def_ms}. Recall that there is an analytic atlas in $M$, and $\bo$ can be assumed to be analytic, too. In other words, in this section, $(M,\bo,g)$ is a real analytic manifold with boundary. 

We will show first that $I_\Gamma f=0$ implies $f^s\in \mathcal{A}(M)$. 
We start with interior analytic regularity. 
Below, $\mathrm{WF}_{\mathrm{A}}(f)$ stands for the analytic wave front set of $f$, see \cite{Sj-Ast,T}.

\begin{proposition}   \label{lemma_wf} 
Let $(x_0,\xi^0)\in T^*M\setminus 0$, and let $\gamma_0$ be a fixed simple geodesic through $x_0$ normal to $\xi^0$. Let  $If(\gamma)=0$ for some 2-tensor $f\in L^2(M)$ and all $\gamma\in \n(\gamma_0)$. 
Let $g$ be analytic in $\n(\gamma_0)$ and $\delta f=0$ near $x_0$. Then  
\be{wf}
(x_0,\xi^0) \not\in\mathrm{WF}_{\mathrm{A}}(f).
\ee
\end{proposition}


\begin{proof}  As explained in Section~\ref{sec_sm}, without loss of generality, we can assume that $\gamma_0$ does not self-intersect. 
Let $U$ be a tubular neighborhood of $\gamma_0$  with $x=(x',x^n)$   analytic semigeodesic coordinates in it, as in the second paragraph of Section~\ref{sec_sm}. 
We can assume that   $x_0=0$, $g_{ij}(0)=\delta_{ij}$, and  $x'=0$ on $\gamma_0$.   In those coordinates, $U$ is given by $|x'|<\eps$, $l^-<x^n<l^+$, with some $0<\eps\ll1$, and we can choose $\eps\ll1$ so that $\{x^n=l^\pm;\; |x'|\le\eps\}$ lie outside $M$. Recall that the lines $x'=\text{const.}$ in $U$ are geodesics.

Then  $\xi^0=((\xi^0)',0)$ with $\xi^0_n=0$. We need to show   that
\be{39}
(0,\xi^0) \not\in \text{WF}_{\text{A}}(f).
\ee
We choose a local chart for the geodesics close to $\gamma_0$. Set first $Z = \{x^n=0;\; |x'|<7\varepsilon/8\}$, and denote the $x'$ variable on $Z$ by $z'$. Then $z'$, $\theta'$ (with $|\theta'|\ll1$) are local coordinates in $\n(\gamma_0)$ determined by $(z',\theta') \to \gamma_{(z',0),(\theta',1)}$.  Each such  geodesic is assumed to be defined on $l^-\le t\le l^+$, the same interval on which $\gamma_0$ is defined. 

Let $\chi_N(z')$, $N=1,2,\dots$,  be a sequence of smooth cut-off functions equal to $1$ for $|z'|\le 3\varepsilon/4$, supported in $Z$, and satisfying the estimates 
\be{N}
\left| \partial^\alpha \chi_N\right|\le (CN)^{|\alpha|}, \quad 
|\alpha|\le N,
\ee
see \cite[Lemma~1.1]{T}.  Set $\theta=(\theta',1)$, $|\theta'|\ll1$,  and multiply 
$$
I f \left(\gamma_{(z',0),\theta}\right) =0
$$
by $\chi_N(z') e^{\mathrm{i}\lambda z'\cdot \xi'}$, where $\lambda>0$, $\xi'$ is in a complex neighborhood of $(\xi^0)'$, and integrate w.r.t.\ $z'$ to get
\be{42}
\iint e^{\lambda \mathrm{i}  z'\cdot \xi'} \chi_N(z')  f_{ij}\left( \gamma_{(z',0),\theta} (t)\right) \dot\gamma_{(z',0),\theta}^i(t) \dot\gamma_{(z',0),\theta}^j (t)\, \d t\, \d z'=0.
\ee

For  $|\theta'|\ll1$, $(z',t)\in Z\times(l^-,l^+)$ are local coordinates near $\gamma_0$ given by $x=\gamma_{(z',0),\theta} (t)$. 

If $\theta'=0$, we have $x=(z',t)$. By a perturbation argument, for $\theta'$ fixed and small enough, $(t,z')$ are analytic local coordinates, depending analytically on $\theta'$.  In particular,  $x=(z'+t\theta',t) + O(|\theta'|)$ but this expansion is not  enough for the analysis below. Performing a change of variables in \r{42}, we get
\be{43}
\int e^{\mathrm{i}\lambda z'(x,\theta')\cdot \xi'}  a_N(x,\theta')   f_{ij}( x) b^i(x,\theta') b^j(x,\theta')\, \d x=0
\ee
for $|\theta'|\ll1$, $\forall\lambda$, $\forall\xi'$, where, for $|\theta'|\ll1$,  the function $(x,\theta') \mapsto a_N$ is analytic and positive for $x$ in a neighborhood of $\gamma_0$, vanishing for  $x\not\in U$, and satisfying \r{N}. The vector field $b$ is analytic on $\supp a_N$, and $b(0,\theta') = \theta$, $a_N(0,\theta')=1$. 

To clarify the arguments that follow, note that if $g$ is Euclidean in $\n(\gamma_0)$, then \r{43} reduces to
$$
\int e^{\mathrm{i}\lambda (\xi',-\theta'\cdot\xi')\cdot x} \chi_N  f_{ij}(x) \theta^i \theta^j\, \d x=0,
$$
where $\chi_N = \chi_N (x'-x^n\theta')$. Then $\xi = (\xi',-\theta'\cdot\xi')$ is perpendicular to $\theta=(\theta',1)$. This implies that 
\be{44}
\int e^{\mathrm{i}\lambda \xi\cdot x} \chi_N   f_{ij}(x) \theta^i (\xi)\theta^j(\xi)\, \d x=0
\ee
for any function $\theta(\xi)$ defined near $\xi^0$, such that $\theta(\xi)\cdot\xi=0$.  This has been noticed and used before if $g$ is close to the  Euclidean metric (with $\chi_N=1$), see e.g., \cite{SU1}. We will assume that $\theta(\xi)$ is analytic. A simple argument (see e.g.\ \cite{Sh,SU1}) shows that a constant symmetric tensor $f_{ij}$ is uniquely determined by  the numbers $f_{ij}\theta^i\theta^j$ for finitely  many $\theta$'s (actually, for $N'=(n+1)n/2$ $\theta$'s); and in any open set on the unit sphere, there are such $\theta$'s.  On the other hand, $f$ is solenoidal near $x_0$. To simplify the argument, assume for a moment that $f$ vanishes on $\bo$ and is solenoidal everywhere. Then $\xi ^i\hat f_{ij}(\xi)=0$. Therefore, combining this with \r{44}, we need to  choose $N=n(n-1)/2$ vectors $\theta(\xi)$, perpendicular to $\xi$, that would uniquely determine the tensor $\hat f$ on the plane perpendicular to $\xi$. To this end, it is enough to know that this choice can be made for $\xi=\xi^0$, then it would be true for $\xi\in \n(\xi^0)$. 
This way, $\xi ^i\hat f_{ij}(\xi)=0$ and the $N$ equations \r{44}  with the so chosen $\theta_p(\xi)$, $p=1,\dots,N$, form a system with a tensor-valued symbol elliptic near $\xi=\xi^0$.  
The $C^\infty$ \PDO\ calculus easily implies the statement of the lemma in the $C^\infty$ category, and the complex stationary phase method below, or the analytic \PDO\ calculus in \cite{T} with appropriate cut-offs in $\xi$, implies the lemma in this special case ($g$ locally Euclidean). 

We proceed with the proof in the general case. Since we will localize eventually near $x_0=0$, where $g$ is close to the Euclidean metric, the special case above serves as a useful guideline. On the other hand, we work near a ``long geodesic'' and the lack of  points conjugate to $x_0=0$ along it will play a decisive role in order to allow us to localize near $x=0$.

Let $\theta(\xi)$ be a vector analytically depending on $\xi$ near $\xi=\xi^0$, such that
\be{th}
\theta(\xi)\cdot\xi=0, \quad \theta^n(\xi)=1, \quad \theta(\xi^0) = e_n.
\ee
Here and below, $e_j$ stand for the vectors $\partial/\partial x^j$. Replace $\theta=(\theta',1)$ in \r{43} by $\theta(\xi)$  (the requirement $|\theta'| \ll1$ is fulfilled for $\xi$ close enough to $\xi^0$), to get
\be{45}
\int e^{\mathrm{i}\lambda \varphi(x,\xi)} \tilde a_N(x,\xi)\tilde  f_{ij}( x)\tilde b^i(x,\xi)\tilde  b^j(x,\xi)\, \d x=0,  \
\ee
where $\tilde a_N$ is analytic  near $\gamma_0\times \{\xi^0\}$, and satisfies \r{N} for $\xi$ close enough to $\xi^0$ and all $x$. Next, $\varphi$, $\tilde b$ are analytic  on $\supp \tilde a_N$ for $\xi$ close  to $\xi^0$. In particular,
$$
\tilde b = \dot\gamma_{(z',0),(\theta'(\xi),1)}(t), \quad t=t(x,\theta'(\xi)), \; z'=z'(x,\theta'(\xi)),
$$
and 
$$
\tilde b(0,\xi) = \theta(\xi), \quad \tilde a_N(0,\xi)=1. 
$$
The phase function is given by 
\be{45a}
\varphi(x,\xi) = z'(x,\theta'(\xi))\cdot \xi'.
\ee
To verify that $\varphi$ is a non-degenerate phase in  $\n(0,\xi^0)$, 
i.e., that $\det \varphi_{x\xi}(0,\xi^0)\not =0$, note first that $z'=x'$ when $x^n=0$, therefore, $(\partial z'/\partial x')(0,\theta(\xi))=\Id$. On the other hand, linearizing near $x^n=0$, we easily get $(\partial z'/\partial x^n)(0,\theta(\xi))=-\theta'(\xi)$. Therefore,
\[
\varphi_x(0,\xi) = (\xi', -\theta'(\xi)\cdot \xi') = \xi
\]
by \r{th}. So we get $\varphi_{x\xi}(0,\xi) = \Id$, which proves the non-degeneracy claim above. In particular, we get that $x\mapsto \varphi_\xi(x,\xi)$ is a local diffeomorphism in $\n(0)$ for $\xi\in\n(\xi^0)$, and therefore injective. We need however a semiglobal version of this along $\gamma_0$ as in the lemma below. For this reason we will make the following special choice of $\theta(\xi)$. Without loss of generality we can assume that 
\[
\xi^0 =e^{n-1}.
\]
Set
\be{45_1}
\theta(\xi) = \bigg( \xi_1,\dots,\xi_{n-2}, -\frac{\xi_1^2+\dots+\xi_{n-2}^2 +\xi_n}{\xi_{n-1}} ,1    \bigg).
\ee
If $n=2$, this reduces to $\theta(\xi) = (-\xi_2/\xi_1,1)$. Clearly, $\theta(\xi)$ satisfies \r{th}.  Moreover, we have
\be{45_2}
\frac{\partial\theta}{\partial \xi_\nu}(\xi^0) = e_\nu, \quad \nu=1,\dots,n-2, \quad \frac{\partial\theta}{\partial \xi_{n-1}}(\xi^0) =0, \quad 
\frac{\partial\theta}{\partial \xi_{n}}(\xi^0) = -e_{n-1},
\ee
In particular, the differential of the map $S^{n-1}\ni \xi \mapsto \theta'(\xi)$ is invertible at $\xi=\xi^0=e^{n-1}$.

\begin{lemma}  \label{lemma_phase} Let $\theta(\xi)$ be as in \r{45_1}, and $\varphi(x,\xi)$ be as in \r{45a}. Then there exists $\delta>0$ such that if 
\[
\varphi_\xi(x,\xi) = \varphi_\xi(y,\xi) 
\]
for  some $x\in U$, $|y|<\delta$, $|\xi-\xi^0|<\delta$, $\xi$ complex,  then $y=x$.
\end{lemma}

\begin{proof}
We will study first the case $y=0$, $\xi=\xi^0$, $x'=0$. Since $\varphi_\xi(0,\xi)=0$, we need to show that $\varphi_\xi((0,x^n),\xi^0)=0$ for $(0,x^n)\in U$ (i.e., for $l^-<x^n< l^+$)  implies  $x^n=0$. 

To compute $\varphi_\xi(x,\xi^0)$, we need first to know $\partial z'(x,\theta')/\partial \theta'$ at $\theta'=0$. Differentiate $\gamma'_{(z',0),(\theta',1)}(t)=x'$ w.r.t.\ $\theta'$, where $t=t(x,\theta')$, $z'=z'(x,\theta')$, to get
\[
\partial_{\theta_\nu} \gamma'_{(z',0),(\theta',1)}(t) + 
\partial_{z'} \gamma'_{(z',0),(\theta',1)}(t)\cdot 
\frac{\partial z'}{\partial \theta_\nu} +\dot \gamma'_{(z',0),(\theta',1)}(t) \frac{\partial t}{\partial \theta_\nu} =0.
\]
Plug $\theta'=0$. Since $\partial t/\partial \theta'=0$ at $\theta'=0$, we get
\[
\frac{\partial z'}{\partial \theta_\nu} = - \partial_{\theta_\nu} \gamma'_{(z',0),(\theta',1)}(x^n)\Big|_{\theta'=0,x'=0} = - J'_\nu(x^n),
\]
where the prime denotes the first $n-1$ components, as usual; $J_\nu(x^n)$ is the Jacobi field along the geodesic $x^n\mapsto \gamma_0(x^n)$ with initial conditions $J_\nu(0)=0$, $DJ_\nu(0)=e_\nu$; and $D$ stands for the covariant derivative along $\gamma_0$. Since $z'((0,x^n),\theta'(\xi^0))=0$, by \r{45a} we then get
\[
\frac{\partial\varphi}{\partial\xi_l} ((0,x^n),\xi^0) = -  \frac{\partial \theta^\mu}{\partial \xi_l}(\xi^0)  J_\mu(x^n)\cdot (\xi^0)'.
\]
By \r{45_2}, (recall that $\xi^0=e^{n-1}$),
\be{45_3}
\frac{\partial\varphi}{\partial\xi_l} ((0,x^n),\xi^0) =
\begin{cases}
-J_l^{n-1}(x^n), & l =1,\dots,n-2,\\
                        0,  & l=n-1,\\
J_{n-1}^{n-1}(x^n),  & l =n,
\end{cases}
\ee
where $J^{n-1}_\nu$ is the $(n-1)$-th component of $J_\nu$. 
Now, assuming that the l.h.s.\ of \r{45_3} vanishes for some fixed $x^n=t_0$, we get that $J_\nu^{n-1}(t_0)=0$, $\nu=1,\dots,n-1$. On the other hand, $J_\nu$ are  orthogonal to $e_n$ because the initial conditions $J_\nu(0)=0$, $DJ_\nu(0)=e_\nu$ are orthogonal to $e_n$, too. Since $g_{in}=\delta_{in}$, this means that $J_\nu^n=0$. Therefore, $J_\nu(t_0)$, $\nu=1,\dots,n-1$, form a linearly dependent system of vectors, thus some non-trivial linear combination $a^\nu J_\nu(t_0)$ vanishes.  Then the solution $J_0(t)$ of the Jacobi equation along $\gamma_0$ with initial conditions $J_0(0)=0$, $DJ_0(0)=a^\nu e_\nu$ satisfies $J(t_0)=0$. Since $DJ_0(0)\not=0$, $J_0$ is not identically zero. Therefore, we get that $x_0=0$ and $x=(0,t_0)$ are conjugate points. Since $\gamma_0$ is a simple geodesic  $x_0$, we must have $t_0=0=x^n$. 

The same proof applies if $x'\not=0$ by shifting the $x'$ coordinates. 

Let now $y$, $\xi$ and $x$ be as in the Lemma. The lemma is clearly true for $x$  in the ball $B(0,\eps_1) = \{|x|<\eps_1\}$, where $\eps_1\ll1$, because $\varphi(0,\xi^0)$ is non-degenerate. On the other hand, $\varphi_\xi(x,\xi)\not= \varphi_\xi(y,\xi)$ for $x\in \bar U\setminus B(0,\eps_1)$,  $y=0$, $\xi=\xi^0$. Hence, we still have $\varphi_\xi(x,\xi)\not= \varphi_\xi(y,\xi)$ for a small perturbation of  $y$ and $\xi$.
\end{proof}

The arguments that follow are close to those in \cite[Section~6]{KSU}. We will apply the complex stationary phase method \cite{Sj-Ast}. For $x$, $y$ as in Lemma~\ref{lemma_phase}, and 
 $|\eta-\xi^0|\le\delta/\tilde C$, $\tilde C\gg2$, $\delta\ll1$, multiply \r{45} by 
$$
\tilde \chi(\xi-\eta)e^{\i \lambda( \i (\xi-\eta)^2/2  -\varphi(y,\xi) )},
$$
where $\tilde \chi$ is the characteristic function of the ball $B(0,\delta)\subset \mathbf{C}^n$,  and integrate w.r.t.\ $\xi$ to get
\be{46aa}
\iint e^{\i\lambda \Phi(y,x,\eta,\xi)}\tilde{\tilde{a}}_N(x,\xi,\eta)  f_{ij}( x) \tilde b^i(x,\xi)  \tilde b^j(x,\xi)\, \d x\, \d \xi=0. 
\ee
Here $\tilde{\tilde{a}}_N = \tilde \chi(\xi-\eta)\tilde a_N$ is another  amplitude, analytic and  elliptic for $x$ close to $0$, $|\xi-\eta|  <\delta/\tilde C$, and
\[
\Phi = -\varphi(y,\xi)+\varphi(x,\xi) +\frac{\i}2 (\xi-\eta)^2.
\]
We study the critical points of $\xi\mapsto \Phi$.  If $y=x$, there is a unique (real) critical point $\xi_{\rm c}=\eta$, and it satisfies  $\Im\Phi_{\xi\xi} >0$ at $\xi= \xi_{\rm c}$. For $y\not=x$, there is no real critical point by Lemma~\ref{lemma_phase}. On the other hand, again by Lemma~\ref{lemma_phase}, there is no (complex) critical point if $|x-y|>\delta/C_1$ with some $C_1>0$, and there is 
a unique complex critical point $\xi_{\rm c}$  if $|x-y|<\delta/C_2$, with some $C_2>C_1$, still non-degenerate if $\delta\ll1$. For any $C_0>0$, if we integrate in \r{46aa} for $|x-y|>\delta/C_0$, and use the fact that $|\Phi_\xi|$ has a positive lower bound (for $\xi$ real), we get
\be{45_5}
\bigg| \iint_{|x-y|>\delta/C_0} e^{\i\lambda \Phi(y,x,\eta,\xi)}\tilde{\tilde{a}}_N(x,\xi,\eta)  f_{ij}( x) \tilde b^i(x,\xi)  \tilde b^j(x,\xi)\, \d x\, \d \xi \bigg| \le  C_3(C_3N/\lambda)^N  +CNe^{-\lambda/C}.
\ee
Estimate \r{45_5} is obtained by integrating $N$ times by parts, using the identity
\[
Le^{\i\lambda \Phi} = e^{\i\lambda \Phi}, \quad L := \frac{\wbar\Phi_\xi\cdot \partial_\xi}{\i\lambda|\Phi_\xi|^2}
\]
as well as using the estimate \r{N}, and the fact that on the boundary of integration in $\xi$, the $e^{\i\lambda\Phi}$ is exponentially small. 
Choose $C_0\gg C_2$. 
Note that  $\Im \Phi>0$ for  $\xi\in\partial (\supp \tilde\chi(\cdot-\eta))$, and $\eta$ as above, as long as $\tilde C\gg1$, and by choosing $C_0\gg1$, we can make sure that $\xi_{\rm c} $ is as close to $\eta$, as we want. 

To estimate \r{46aa} for $|x-y|<\delta/C_0$, set 

$$
\psi(x,y,\eta) := \Phi\big|_{\xi=\xi_{\text{c}}}.
$$
Note that $\xi_{\text{c}} =-\i(y-x)+\eta+O(\delta)$, and $\psi(x,y,\eta) = \eta\cdot (x-y) +\frac{\i}2 |x-y|^2+O(\delta)$. We will not use this to study the properties of $\psi$, however. Instead, observe that at $y=x$ we have
\be{46a}
\psi_y(x,x,\eta) = -\varphi_x(x,\eta), \quad
\psi_x(x,x,\eta) =  \varphi_x(x,\eta), \quad
\psi(x,x,\eta)=0.
\ee
We also get that 
\be{46b}
\Im \psi(y,x,\eta) \ge |x-y|^2/C.
\ee
The latter can be obtained by setting $h=y-x$ and expanding in powers of $h$. 
The stationary complex phase method \cite{Sj-Ast}, see Theorem~2.8 there and the remark after it,  gives
\be{47}
\int_{|x-y|\le \delta/C_0} e^{\i \lambda \psi(x,\alpha)}  
  f_{ij}( x) B^{ij}(x,\alpha; \lambda)  \, \d x =  O\big( \lambda^{n/2}(C_3N/\lambda)^N +Ne^{-\lambda/C} \big),  \quad \forall N,
\ee 
where $\alpha = (y,\eta)$, and $B$ is a classical analytic symbol \cite{Sj-Ast} with principal part equal to $\tilde b\otimes \tilde b$, up to an elliptic factor. The l.h.s.\ above is independent of $N$, and choosing $N$ so that $N\le \lambda/(C_3e)\le N+1$ to conclude that the r.h.s.\ above is $O(e^{-\lambda/C})$.

In preparation for applying the characterization of an analytic wave front set through a generalized FBI transform \cite{Sj-Ast}, define the transform
$$
\alpha \longmapsto \beta = \left(\alpha_x, \nabla_{\alpha_x}\varphi(\alpha)\right),
$$
where, following \cite{Sj-Ast}, $\alpha=(\alpha_x,\alpha_\xi)$. It is a diffeomorphism from $\n(0,\xi^0)$ to its image, and denote the inverse one by $\alpha(\beta)$. Note that this map and its inverse  preserve the first (n-dimensional) component and change only the second one.  This is equivalent to setting $\alpha=(y,\eta)$, $\beta = (y,\zeta)$, where  $\zeta = \varphi_y(y,\eta)$. Note that $\zeta =\eta+O(\delta)$, and at $y=0$, we have $\zeta=\eta$.

Plug $\alpha=\alpha(\beta)$ in \r{47} to get 
\be{48}
\int e^{\i \lambda \psi(x,\beta)} 
 f_{ij}( x)B^{ij}(x,\beta; \lambda) \, \d x =  O\big( e^{-\lambda/C} \big),  
\ee
where $\psi$,  $B$ are (different) functions having the same properties as above. Then 
\be{49}
\psi_y(x,x,\zeta) = -\zeta, \quad
\psi_x(x,x,\zeta) =  \zeta, \quad
\psi(x,x,\zeta)=0.
\ee
The symbols in \r{48} satisfy
\be{48_1}
\sigma_p(B)(0,0,\zeta) \equiv\theta(\zeta)\otimes \theta(\zeta) \quad \text{up to an elliptic factor},
\ee
and in particular, $\sigma_p(B)(0,0,\xi^0)\equiv e_n\otimes e_n$, where $\sigma_p$ stands for the principal symbol. 

Let $\theta_1=e_n, \, \theta_2, \dots,\theta_N$ be $N=n(n-1)/2$ unit vectors at $x_0=0$, normal to $\xi^0=e^{n-1}$ such that any constant symmetric 2-tensor $f$ such that $f_i^{n-1}=0$, $\forall i$ (i.e., $f_i^j\xi^0_j=0$) is uniquely determined by  $f_{ij}\theta^i\theta^j$, $\theta=\theta_p$, $p=1,\dots,N$. Existence of such vectors is easy to establish, as mentioned above, and one can also see that such a set exists in any open set in $(\xi^0)^\perp$. We can therefore assume that $\theta_p$ belong to a small enough neighborhood of $\theta_1=e_n$ such that the geodesics $[-l^-,l^+] \ni t\mapsto \gamma_{0,\theta_p}(t)$  through $x_0=0$ are all simple. Then we can rotate a bit the coordinate system such that $\xi^0=e^{n-1}$ again, and $\theta_p=e_n$, and repeat the construction above. This gives us $N$ phase functions $\psi_{(p)}$, and as many symbols  $B_{(p)}$ in \r{48} such that \r{49} holds for all of them, i.e., in the coordinate system related to $\theta_1=e_n$, we have
\be{48_1a}
\int e^{\i \lambda \psi_{(p)}(x,\beta)}  
 f_{ij}( x)B^{ij}_{(p)}(x,\beta; \lambda)\, \d x =  O\big( e^{-\lambda/C} \big),  \quad p=1,\dots,N,
\ee
and by \r{48_1},
\be{48_2}
\sigma_p(B_{(p)})(0,0,\xi^0) \equiv \theta_p \otimes\theta_p, \quad p=1,\dots,N.
\ee

Recall that   $\delta f=0$ near $x_0=0$. Let $\chi_0=\chi_0(x)$ be a smooth cutoff close enough to $x=0$, equal to $1$ in $\n(0)$. Integrate $\frac1{\lambda} \exp\big(\i \lambda \psi_{(1)}(x,\beta)\big) \chi_0\delta f =0$ w.r.t.\ $x$, and by \r{46b}, after an integration by parts, we get
\be{49a}
\int e^{\i \lambda \psi_{(1)}(x,\beta)} \chi_0(x)f_{ij}(x)  C^j(x,\beta;\lambda)\, \d x= O\big( e^{-\lambda/C} \big), \quad i=1,\dots,n,
\ee
for $\beta_x=y$ small enough, where $\sigma_p(C^j)(0,0,\xi^0)=(\xi^0)^j$. 

Now, the system of $N+n = (n+1)n/2$ equations \r{48_1a}, \r{49a} can be viewed as a tensor-valued operator applied to the tensor $f$. Its symbol, an elliptic factor at $(0,0,\xi^0)$, has ``rows'' given by $\theta_p^i \theta_p^j$, $p=1,\dots,N$; and $\delta^i_k(\xi^0)^j$, $k=1,\dots,n$. It is easy to see that it is elliptic; indeed, the latter is equivalent to the statement that if for  some (constant) symmetric 2-tensor $f$, in Euclidean geometry (because $g_{ij}(0)=\delta_{ij}$), we have $f_{ij}\theta_p^i \theta_p^j=0$, $p=1,\dots,N$; and $f_i^{n-1}=0$, $i=1,\dots,n$, then $f=0$. This however follows from the way we chose $\theta_p$. Therefore, \r{39} is a consequence of \r{48_1a}, \r{49a}, see  \cite[Definition~6.1]{Sj-Ast}. Note that in \cite{Sj-Ast}, it is  required that $f$ must be replaced by $\bar f$ in \r{48_1a}, \r{49a}. If $f$ is complex-valued, we could use the fact that $I(\Re f)(\gamma)=0$, and $I(\Im f)(\gamma)=0$ for $\gamma$ near $\gamma_0$ and then work with real-valued $f$'s only. 

Since the phase functions in \r{48_1a} depend on $p$, we need to explain why the characterization of the analytic wave front sets in \cite{Sj-Ast} can be generalized to this vector-valued case. The needed modifications are as follows. We define 
$h^{ij}_{(p)}(x,\beta;\lambda) = B_{(p)}^{ij}$, $p=1,\dots,N$; and $h^{ij}_{(N+k)}(x,\beta;\lambda) = C^{j}\delta^i_{k}$, $k=1,\dots,n$. Then $\{h^{ij}_{(p)}\}$, $p=1,\dots,N+n$, is an elliptic symbol near $(0,0,\xi^0)$. In the proof of \cite[Prop.~6.2]{Sj-Ast}, under the conditions \r{46b}, \r{49}, the operator $Q$ given by 
$$
[Qf]_p(x,\lambda) = \iint e^{\i \lambda( \psi_{(p)}(x,\beta) - \overline{\psi_{(p)}(y,\beta)}  )} 
 f_{ij}( y,\lambda) h_{(p)}^{ij}(x,\beta; \lambda) \, \d y\, \d \beta 
$$
is a \PDO\ in the complex domain with an elliptic matrix-valued symbol, where we view $f$ and $Qf$ as vectors in $\R^{N+n}$. Therefore, it admits a parametrix in $H_{\psi,x_0}$ with a suitable $\psi$ (see \cite{Sj-Ast}). Hence, one can find an analytic classical matrix-valued symbol $r(x,\beta,\lambda)$ defined near $(0,0,\xi^0)$, such that for any constant symmetric $f$ we have
$$
\left[Q\left( r(\cdot,\beta,\lambda) e^{\i \lambda\psi_{(1)}}  f\right) \right]_p =
e^{\i \lambda\psi_{(1)}}f, \quad \forall p.
$$
The rest of the proof is identical to that of \cite[Prop.~6.2]{Sj-Ast} and allows us to show that \r{48} is preserved with a different choice of the phase functions satisfying \r{46b}, \r{49}, and elliptic amplitudes; in particular,
$$
\int e^{\i \lambda\psi_{(1)}(x,\beta)} \chi_2(x)  f_{ij}(x) \, \d x = O\big( e^{-\lambda/C} \big) , \quad \forall i,j
$$
{}for $\beta\in\n(0,\xi^0)$ and for some standard cut-off $\chi_2$ near $x=0$. This proves \r{39}, see \cite[Definition~6.1]{Sj-Ast}.

This concludes the proof of Proposition~\ref{lemma_wf}. Notice that the proof works in the sane way, if $f$ is a distribution valued tensor field, supported in $M$. 
\end{proof}

\begin{lemma} \label{pr_an}
Under the assumptions of Theorem~\ref{thm_an}, let $f$ be such that $I_\Gamma f=0$.
Then $f^s\in\mathcal{A}(M)$. 
\end{lemma}

\begin{proof} Proposition~\ref{lemma_wf}, combined with the completeness of $\Gamma$, imply that $f^s$ is analytic in the interior of $M$. To prove analyticity up to the boundary, we do the following.

We can assume that $M_1\setminus M$ is defined by $-\eps_1\le x^n\le 0$, where $x^n$ is a boundary normal coordinate. Define  the manifold $M_{1/2}\supset M$ by $x^n\ge -\eps_1/2$, more precisely, $M_{1/2} = M\cup \{-\eps_1/2\le x^n\le0\}\subset M_1$. 

We will show first that $f^s_{M_{1/2}}\in \mathcal{A}(M_{1/2})$. Let us first notice, that in $M_{1/2}\setminus M$, $f^s_{M_{1/2}} = -dv_{M_{1/2}}$, where $v_{M_{1/2}}$ satisfies $\Delta^s v_{M_{1/2}}=0$ in $M_{1/2}\setminus M$, $v|_{\bo_{1/2}}=0$. Therefore, $v_{M_{1/2}}$ is analytic up to $\partial M_{1/2}$ in   $M_{1/2}\setminus M$, see \cite{MN, SU-rig}. Therefore, we only need to show that $f^s_{M_{1/2}}$ is analytic in some neighborhood 
of $M$. This however follows from Proposition~\ref{lemma_wf}, applied to $M_{1/2}$. Note that if $\eps_1\ll1$, simple geodesics through some $x\in M$ would have endpoints outside $M_{1/2}$ as well, and by a compactness argument, we need finitely many such geodesics to show that Proposition~\ref{lemma_wf} implies that $f^s_{M_{1/2}}$ is analytic in, say, $M_{1/4}$, where the latter is defined similarly to $M_{1/2}$ by $x^n\ge -\eps_1/4$.  

To compare $f^s_{M_{1/2}}$ and $f^s = f^s_M$, see also \cite{SU-Duke, SU-rig}, write $f^s_{M_{1/2}} = f-dv_{M_{1/2}}$ in $M_{1/2}$, and $f^s_M = f-dv_M$ in $M$. Then $dv_{M_{1/2}}= -f^s_{M_{1/2}}$ in $M_{1/2}\setminus M$, and is therefore analytic there, up to  $\bo$. Given $x\in\bo$, integrate $\langle dv_{M_{1/2}}, \dot\gamma^2 \rangle$ along geodesics in $M_{1/2}\setminus M$, close to ones normal to the boundary, with initial point $x$ and endpoints on $\bo_{1/2}$. Then we get that $v_{M_{1/2}}|_{\bo} \in \mathcal{A}(\bo)$. Note that $v_{M_{1/2}}\in H^1$ near $\bo$, and taking the trace on $\bo$ is well defined, and moreover, if $x^n$ is a boundary normal coordinate, then $\n(0)\ni x^n \mapsto v_{M_{1/2}}(\cdot,x^n)$ is  continuous. Now,
\be{50}
f^s_M = f-dv_M = f^s_{M_{1/2}} +dw \quad \mbox{in $M$,}\quad \mbox{where $w = v_{M_{1/2}}-v_M$.}
\ee
The vector field $w$ solves
$$
\upDelta^s w = 0, \quad w|_{\bo} = v_{M_{1/2}}|_{\bo} \in\mathcal{A}(\bo).
$$
Therefore, $w\in \mathcal{A}(M)$, and by \r{50}, $f^s_M\in \mathcal{A}(M)$. 

This completes the proof of Lemma~\ref{pr_an}.
\end{proof}

\begin{proof}[Proof of Theorem~\ref{thm_an}] Let $I_\Gamma f=0$. 
We can assume first that $f=f^s$, and then $f\in\mathcal{A}(M)$ by Lemma~\ref{pr_an}. 
By Lemma~\ref{lemma_bd}, there exists $h \in \mathcal{S}^{-1} \mathcal{S} f$ such that $\partial^\alpha h=0$ on $\bo$ for all $\alpha$. The tensor field $h$ satisfies \r{l1_2}, i.e., $h_{ni}=0$, $\forall i$, in boundary normal coordinates, which is achieved by setting $h=f-dv_0$, where $v_0$ solves \r{a1} near $\bo$. Then $v_0$, and therefore, $h$ is analytic for small $x^n\ge0$, up to $x^n=0$. Lemma~\ref{lemma_bd} then implies that $h=0$ in $\n(\bo)$. 
So we get that
\be{50a}
f=dv_0, \quad 0\le x^n<\eps_0,\quad \text{with $ v_0|_{x^n=0}=0$},
\ee
where $x^n$ is a global normal coordinate, and $0<\eps_0\ll1$. Note that the solution $v_0$ to \r{50a} (if exists, and in this case we know it does) is unique, as can be easily seen by integrating $\langle f,\dot\gamma^2\rangle$ along paths close to normal ones to $\bo$ and using \r{v}.

We show next that $v_0$ admits an analytic continuation from a neighborhood of any $x_1\in\bo$ along any path in $M$. 

Fix $x\in M$. Let $c(t)$, $0\le t\le1$ be a path in $M$ such that $c(0)=x_0 \in\bo$ and $c(1)=x$. Given $\eps>0$, one can find a polygon $x_0x_1\dots x_k x$ consisting of geodesic segments of length not exceeding $\eps$, that is close enough and therefore homotopic to $c$. One can also assume that the first one is transversal to $\bo$, and if $x\in\bo$, the last one is transversal to $\bo$ as well; and all other points of the polygon are in $\Mint$. We choose $\eps\ll1$ so that there are no conjugate points on each geodesic segment above. We also assume that $\eps\le\eps_0$. Then $f=dv$  near $x_0x_1$  with $v=v_0$ by \r{50a}. As in the second paragraph of Section~\ref{sec_sm}, one  can choose  semigeodesic coordinates $(x',x^n)$  near $x_1x_2$, and a small enough hypersurface $H_1$ through $x_1$  given locally by $x^n=0$. 
As in Lemma~\ref{lemma_bd},  one can find an analytic 1-form $v_1$ defined near $x_1x_2$, so that $(f-dv_1)_{in}=0$, $v_1|_{x^n=0}=v_0(x',0)$.  Close enough to $x_1$, we have $v_1=v_0$ because $v_0$ is also a solution, and the solution is unique, see also \r{a1'}. 
Since $v_1$ is analytic, we get that it is an analytic extension of $v_0$ along $x_1x_2$. Since $f$ and $v_1$ are both analytic in $\n(x_1x_2)$, and $f=dv_1$ near $x_1$, this is also true in $\n(x_1x_2)$. So we extended $v_0$ along $x_0x_1x_2$, let us call this extension $v$. 
Then we do the same thing near $x_2x_3$, etc., until we reach $\n(x)$, and then $f=dv$ there.

This defines $v$ in $\n(x)$, where $x\in M$ was chosen arbitrary. It remains to show that this definition is independent of the choice of the path. Choose another path that connects some $y_1\in\bo$ and $x$.  Combine them both to get a path that connects $x_1\in \bo$ and $y_1\in\bo$. It suffices to prove that the analytic continuation of $v_0$ from $x_1$ to $y_1$ equals $v_0$ again.  Let $c_1\cup \gamma_1 \cup c_2\cup\gamma_2\cup\dots\cup  \gamma_k \cup c_{k+1}$ be the polygon homotopic to the path above. Analytic continuation along  $c_1$ coincides with $v_0$ again by \r{50a}. Next, let $p_1$, $p_2$ be the initial and the endpoint of $\gamma_1$, respectively, where $p_1$ is also the endpoint of $c_1$. 
We continue analytically $v_0$ from $\n(p_1)$ to $\n(p_2)$ along $\gamma_1$, let us call this continuation $v$. By what we showed above, $f=dv$ near $\gamma_1$. Since $If(\gamma_1)=0$, and $v(p_1)=0$, we get by \r{v}, that $\langle v(p_2),\dot \gamma_1(l)\rangle =0$ as well, where $l$ is such $\gamma_1(l)=p_2$. Using the assumption that $\gamma_1$ is transversal to $\bo$ at both ends, one can perturb the tangent vector $\dot \gamma_1(l)$ and this will define a new geodesic through $p_2$ that hits $\bo$ transversely again near $p_1$, where $v=v_0=0$.  Since $\Gamma$ is open, integral of $f$ over this geodesic vanishes again, therefore $\langle v(p_2),\xi_2\rangle =0$ for $\xi_2$ in an open set. Hence $v(p_2)=0$. Choose $q_2\in\bo$ close enough to $p_2$, and $\eta_2$ close enough to $\xi_2$ (in a fixed chart). Then the geodesic through $(q_2,\eta_2)$ will hit $\bo$ transversally close to $p_1$, and we can repeat the same arguments. We therefore showed that $v=0$ on $\bo$ near $p_2$. On the other hand, $v_0$ has the same property. Since $f=dv=dv_0$ there, by the remark after \r{50a}, we get that $v=v_0$ near $p_2$. We repeat this along all the legs of the polygon until we get that the analytic continuation $v$ of $v_0$ along the polygon, from $x_1$ to $y_1$, equals $v_0$ again. 

As a consequence of this, we get that $f=dv$ in $M$ with $v=0$ on $\bo$. Since $f=f^s$, this implies $f=0$.

This completes the proof of Theorem~\ref{thm_an}.
\end{proof}

\section{Proof of Theorems~\ref{thm_stab} and \ref{thm_I}}
\begin{proof}[Proof of Theorem~\ref{thm_stab}]
Theorem~\ref{thm_stab}(b), that also implies (a), is a consequence of Proposition~\ref{pr_2}, as shown in \cite{SU-rig}, see the proof of Theorem~2 and Proposition~4 there. Part (a) only follows more directly from \cite[Prop.~V.3.1]{Ta1} and its generalization, see \cite[Thm~2]{SU-Duke}.
\end{proof}

\begin{proof}[Proof of Theorem~\ref{thm_I}]
First, note that for any analytic metric in $\mathcal{G}$, $I_{\Gamma_g}$ is s-injective  by Theorem~\ref{thm_an}. 
We  build $\mathcal{G}_s$ as a small enough neighborhood of the analytic metrics in $\mathcal{G}$. 
 Then $\mathcal{G}_s$ is dense in $\mathcal{G}$ (in the $C^k(M_1)$ topology) since it includes the analytic metrics. To complete the definition of $\mathcal{G}_s$, fix an analytic $g_0\in \mathcal{G}$. By Lemma~\ref{lemma_H}, one can find  $\mathcal{H}'\Subset\mathcal{H}$ related to $g=g_0$ and $\Gamma_g$,  satisfying the assumptions of Theorem~\ref{thm_stab}, and they have the properties required for $g$ close enough to $g_0$.

Let $\alpha$ be as in Theorem~\ref{thm_stab} with $\alpha=1$ on $\mathcal{H}'$. Then, by Theorem~\ref{thm_stab}, $I_{\alpha,g}$ is s-injective for $g$ close enough to $g_0$ in $C^k(M_1)$. By Lemma~\ref{lemma_1}, for any such $g$, $I_{\Gamma^\alpha}$ is s-injective, where $\Gamma^\alpha = \Gamma(\mathcal{H}^\alpha)$, $\mathcal{H}^\alpha = \supp\alpha$. If $g$ is close enough to $g_0$, $\Gamma^\alpha\subset \Gamma_g$ because when $g=g_0$, $\Gamma^\alpha\subset \Gamma(\mathcal{H})\Subset\Gamma_{g_0}$, and $\Gamma_g$ depends continuously on $g$ in the sense described before the formulation of Theorem~\ref{thm_I}. Those arguments show that there is a neighborhood of each analytic $g_0\in\mathcal{G}$ with an s-injective $I_{\Gamma_g}$. Therefore, one can choose an open dense subset $\mathcal{G}_s$ of $\mathcal{G}$ with the same property.
\end{proof}

\begin{proof}[Proof of Corollary~\ref{cor_1}.]
It is enough to notice that the set of all simple geodesics  related to $g$ 
depends conti\-nuously on $g$ in the sense of Theorem~\ref{thm_I}. Then the proof follows from the paragraph above.
\end{proof}

\section{X-ray transform of functions and 1-forms/vector fields}  \label{sec_f}

If $f$ is a vector field on $M$, that we identify with an  1-form, then its X-ray transform is defined quite similarly to \r{I_G} by 
\be{I_v}
I_\Gamma f(\gamma) = \int_0^{l_\gamma} \langle f(\gamma(t)), \dot \gamma(t) \rangle \,\d t,  \quad \gamma\in \Gamma.
\ee
If $f$ is a function on $M$, then we set
\be{I_f}
I_\Gamma f(\gamma) = \int_0^{l_\gamma} f(\gamma(t))\,\d t,  \quad \gamma\in \Gamma.
\ee
The latter case is a partial case of the X-ray transform of 2-tensors; indeed, if $f =\alpha g$, where $f$ is a 2-tensor, $\alpha$ is a function, and $g$ is the metric, then $I_\Gamma f = I_\Gamma\alpha$, where in the l.h.s., $I_\Gamma$ is as in \r{I_G}, and on the right, $I_\Gamma$ is as in \r{I_f}. The proofs for the X-ray transform of functions are  simpler, however, and in particular, there is no loss of derivatives in the estimate \r{est}, as in \cite{SU-Duke}. This is also true for the X-ray transform of vector fields and the proofs are more transparent than those for tensors of order 2 (or higher). Without going into details (see \cite{SU-Duke} for the case of simple manifolds), we note that  the main theorems in the Introduction remain true. In case of 1-forms, estimate \r{est} can be improved to
\be{est1}
\|f^s\|_{L^2(M)}/C  \le 
\|N_{\alpha} f\|_{H^1(M_1)} \le C\|f^s\|_{L^2(M)},
\ee
while in case of functions, we have 
\be{est2}
\|f\|_{L^2(M)}/C  \le 
\|N_{\alpha} f\|_{H^1(M_1)} \le C\|f\|_{L^2(M)}.
\ee
If $(M,\bo)$ is simple, then the full X-ray transform of functions and 1-forms (over all geodesics) is injective, respectively s-injective, see \cite{Mu2, Mu-R, BG, AR}.

\end{document}